\documentclass[preprint,12pt]{elsarticle}

\usepackage{amssymb}
\usepackage{pstricks}
\usepackage{amssymb}
\usepackage{float,afterpage}
\usepackage{palatino}
\usepackage{leftidx}
\usepackage{srcltx}
\usepackage{soul}
\usepackage{mathptmx}
\usepackage{ amsthm}
\usepackage{amsmath}
\usepackage{pstricks-add}
\usepackage{bm}
\usepackage{color}
\usepackage{amsfonts}
\usepackage{cancel}
\newtheorem{thm}{Theorem}
\newtheorem{defn}[thm]{Lemma}
\newtheorem{prop}[thm]{Lemma}
\newtheorem{lem}[thm]{Lemma}
\newtheorem{cor}[thm]{Lemma}
\newtheorem{rem}[thm]{Remark}
\def\bbR{{\mathbb R}}

\def\bbN{{\mathbb N}}

\def\a{{\~{a}}}

\begin{document}

\begin{frontmatter}



\title{Scaling limit for a family of random paths with radial behavior}

\author[label1,label2]{Cristian F. Coletti, Leon A. Valencia}
\address[label1]{Centro de Matem\'atica, Computa\c c\a o e Cogni\c c\a o - Universidade Federal do ABC. Avenida do Estado 5001, Santo Andr\'e, S\a o Paulo, Brasil.}
\address[label2]{Instituto de Matem\'aticas - Universidad de Antioquia. Calle 67 Nº 53-108, Medell\'in, Colombia.}

%\author{fgds}

%\address{fgf}

\begin{abstract}
We introduce a system of coalescing random paths with radial behavior in a subset of the plane. We call it the {\it{Discrete Radial Poissonian Web}}. 
We show that under diffusive scaling this family converges in distribution to a mapping of a restriction of the Brownian Web.
\end{abstract}


% \begin{highlights}
% \end{highlights}


\end{frontmatter}


% \SECTION{}
% \LABEL{}



%%       \citet{<label>} ==> Jones et al. [21]
%%       \citep{<label>} ==> [21]
%%



\section{Introduction}
%%%%%%%%%%%%%%%%%%%%%%%%%%%%%%%%%%%%%%%%%%%%%%%%%%%%%%%%%%%%%%%
\label{sec:intro}
The aim of this work is to introduce a family of coalescing random paths with radial behavior and study its scaling limit in the diffusive scaling. The scaling limit of
the Radial Web considered in this paper is related to the Brownian Web. The Brownian Web (see \cite{FINR}) is, roughly speaking, a family of coalescing Brownian motions
starting from every point in $\mathbb{R}^2$. Indeed, the scaling limit of our model belongs the family of the {\it Brownian Bridge Web } introduced in \cite{valle2015}. In that work the authors considered a variant of the radial spanning tree introduced in \cite{baccelli2007} which is much closer to the radial spanning tree than our model. They showed that locally the rescaled family of random paths converges in distribution to the
%they obtained its scaling limit (in the diffusive scaling) and called it the 
{\it Brownian Bridge Web}. Despite our model being simpler, it also converges in distribution to the Brownian Bridge Web under the diffusive scaling. Thus we provide a new family of coalescing random paths converging to the Brownian Bridge Web. It is worth mentioning that the simplicity of our model allows us to use, in a smart way, the FKG inequality to show condition $B_2$ (see section \ref{vcb2}).

We begin with a tentative description of the Discrete Radial Web. Consider $n$ circles centered at the origin with radius $n,n-1, ...,1$ respectively. Define the next
circle of a circle centered at the origin with radius $k$ as the circle centered at the origin with radius $k-1$.
On each of these circles, consider a Poisson point process of rate $1$. Assume independence between these point processes. Then, from each Poissonian point, draw a line
to the nearest Poissonian point on the next circle, if any. If such a point does not exist then, connect it to the nearest point on the next circle to the previous one.
Repeat this procedure until connecting all paths to the origin. Proceeding in this way we obtain a family of {\it{coalescing random paths}}. See figure~{\ref{fig:figure1}} 
in order to figure out this description. This family is what we would like to call the Discrete
Radial Poissonian Web (DRPW). For technical reasons we need to make some restrictions. In next section we provide a formal description of the mapping of the restricted Brownian Web 
which we call the $T$ - Brownian Web and state
the main result of this work, namely the weak convergence of the DRPW in the diffusive scaling. To prove the convergence result we introduce and verify a convergence criteria which is
entirely analogous to the convergence criteria given in \cite{FINR}. Following the notation introduced in that work, we call the conditions of the convergence
criteria $I, B_1$ and $B_2$. This criteria has been verified in a number of papers, for instance, see \cite{FFW}, \cite{SUN} and \cite{CFD} . In \cite{CFD}, the authors proved 
convergence to the Brownian Web for a system of random walks introduced in \cite{GRS}. In \cite{CV}, the authors proposed a generalization of this system of random walks by 
considering a system of coalescing nonsimple random walks and proved convergence to the Brownian web. The main difference between the system of random paths considered in \cite{CV} and 
the one considered in the present work is due to the radial behavior of its paths.

This paper is organized as follows. In section $2$ we present the basic definitions, notations and state the main result of this work. The end of section $2$ is devoted
to the study of the DRPW and to the proof of the main result of this work. In section $3$ we state and prove the result about the asymptotic behavior of the tail of the
coalescence time for two random paths of the DRPW. In section $4$ we verify condition $B_2$. Section $5$ is entirely devoted to the verification of condition $I$.
Condition $B_1$ is verified in section $6$. For the sake of completeness we include, at the end of this work, an appendix about weak convergence.

\vspace{1cm}
\definecolor{xdxdff}{rgb}{0.0, 0.5, 1.0}.
\psset{unit=0.6}
\begin{figure}
\label{figure1}
\begin{pspicture*}(-3.94,-6.56)(16.2,5.88)
\pscircle(6.02,-0.14){6}
\pscircle(6.02,-0.14){5.5}
\pscircle(6.02,-0.14){5}
\pscircle(6.02,-0.14){4.5}
\pscircle(6.02,-0.14){4}
\pscircle(6.02,-0.14){3.5}
\pscircle(6.02,-0.14){3}
\pscircle(6.02,-0.14){2.5}
\pscircle(6.02,-0.14){2}
\pscircle(6.02,-0.14){1.5}
\pscircle(6.02,-0.14){1}
\pscircle(6.02,-0.14){0.5}
\rput[tl](3.28,-5.1){\color{blue}}
\rput[tl](6.68,-5.72){\color{blue}}
\rput[tl](4.58,-5.82){\color{blue}}
\rput[tl](11.24,-2.88){\color{blue}}
\rput[tl](9.88,-4.54){\color{blue}}
\rput[tl](1.1,-3.54){\color{blue}}
\rput[tl](0.06,3.36){\color{blue}}
\rput[tl](0.1,1.28){\color{blue}}
\rput[tl](-0.02,-0.48){\color{blue}}
\rput[tl](8.6,5.4){\color{blue}}
\rput[tl](3.54,5.5){\color{blue}}
\rput[tl](10.6,3.86){\color{blue}}
\rput[tl](11.5,2.38){\color{blue}}
\rput[tl](11.94,-0.48){\color{blue}}
\rput[tl](0.62,-1.22){\color{blue}}
\rput[tl](1.04,2.44){\color{blue}}
\rput[tl](2.32,4.14){\color{blue}}
\rput[tl](3.2,-4.74){\color{blue}}
\rput[tl](1.18,-2.72){\color{blue}}
\rput[tl](10.92,-2.4){\color{blue}}
\rput[tl](4.48,5.32){\color{blue}}
\rput[tl](9.52,-4.2){\color{blue}}
\rput[tl](6.16,-5.5){\color{blue}}
\rput[tl](7.32,5.32){\color{blue}}
\rput[tl](3.64,-4.42){\color{blue}}
\rput[tl](5.44,4.98){\color{blue}}
\rput[tl](10.78,1.44){\color{blue}}
\rput[tl](7.22,-4.8){\color{blue}}
\rput[tl](10.86,-0.88){\color{blue}}
\rput[tl](2.1,3.14){\color{blue}}
\rput[tl](8.4,4.38){\color{blue}}
\rput[tl](4.06,4.62){\color{blue}}
\rput[tl](9.94,-3.04){\color{blue}}
\rput[tl](10.26,-1.3){\color{blue}}
\rput[tl](1.42,-2.08){\color{blue}}
\rput[tl](7.96,-4.02){\color{blue}}
\rput[tl](9.16,3.16){\color{blue}}
\rput[tl](1.62,1.22){\color{blue}}
\rput[tl](0.08,0.4){\color{blue}}
\rput[tl](6.62,4.46){\color{blue}}
\rput[tl](5.02,-4.4){\color{blue}}
\rput[tl](2.18,-2.4){\color{blue}}
\rput[tl](2.52,2.9){\color{blue}}
\rput[tl](2.74,2.36){\color{blue}}
\rput[tl](3.86,3.4){\color{blue}}
\rput[tl](7.08,3.86){\color{blue}}
\rput[tl](9.66,-1.52){\color{blue}}
\rput[tl](8.4,3.18){\color{blue}}
\rput[tl](3.88,-3.42){\color{blue}}
\rput[tl](9.28,2.22){\color{blue}}
\rput[tl](6.72,-3.92){\color{blue}}
\rput[tl](5.8,-3.5){\color{blue}}
\rput[tl](1.98,-0.4){\color{blue}}
\rput[tl](6.98,3.36){\color{blue}}
\rput[tl](8.04,2.16){\color{blue}}
\rput[tl](9.08,-1.58){\color{blue}}
\rput[tl](3,0.42){\color{blue}}
\rput[tl](8.96,0.08){\color{blue}}
\rput[tl](2.66,-1.14){\color{blue}}
\rput[tl](2.68,1.2){\color{blue}}
\rput[tl](8.98,1.76){\color{blue}}
\rput[tl](6.56,-2.94){\color{blue}}
\rput[tl](7.58,-1.18){\color{blue}}
\rput[tl](3.78,-1.22){\color{blue}}
\rput[tl](3.96,2.24){\color{blue}}
\rput[tl](6.6,2.42){\color{blue}}
\rput[tl](6.26,0.38){\color{blue}}
\rput[tl](3.96,0.1){\color{blue}}
\psline(3.33,-5.5)(3.27,-4.9)
\psline(3.27,-4.9)(3.69,-4.56)
\psline(3.69,-4.56)(5.08,-4.54)
\psline(5.08,-4.54)(3.93,-3.55)
\psline(3.93,-3.55)(5.86,-3.64)
\psline(5.86,-3.64)(6.62,-3.08)
\psline(6.62,-3.08)(3.83,-1.34)
\psline(3.83,-1.34)(4.02,0)
\psline(4.02,0)(6.32,0.26)
\psline(6.32,0.26)(6.02,-0.14)
\psline(4.7,-5.97)(6.22,-5.64)
\psline(6.22,-5.64)(7.3,-4.97)
\psline(7.3,-4.97)(8.02,-4.17)
\psline(8.02,-4.17)(6.78,-4.07)
\psline(6.78,-4.07)(5.86,-3.64)
\psline(7.75,-5.89)(6.22,-5.64)
\psline(9.94,-4.68)(9.57,-4.34)
\psline(9.57,-4.34)(9.98,-3.2)
\psline(9.98,-3.2)(10.33,-1.42)
\psline(10.33,-1.42)(9.72,-1.65)
\psline(9.72,-1.65)(9.13,-1.75)
\psline(9.13,-1.75)(9.02,0)
\psline(9.02,0)(6.64,2.28)
\psline(6.64,2.28)(4.02,0)
\psline(7.63,-1.33)(6.32,0.26)
\psline(11.28,-3.03)(10.96,-2.57)
\psline(10.96,-2.57)(9.98,-3.2)
\psline(10.95,-0.99)(10.33,-1.42)
\psline(12,-0.68)(10.96,-2.57)
\psline(11.53,2.24)(8.48,4.21)
\psline(8.48,4.21)(9.21,3.03)
\psline(9.21,3.03)(8.44,3.04)
\psline(8.44,3.04)(7.05,3.21)
\psline(7.05,3.21)(8.06,2.06)
\psline(8.06,2.06)(6.64,2.28)
\psline(10.81,1.29)(9.21,3.03)
\psline(9.35,2.08)(7.05,3.21)
\psline(9.04,1.63)(8.06,2.06)
\psline(10.63,3.7)(7.34,5.2)
\psline(7.34,5.2)(5.5,4.83)
\psline(5.5,4.83)(6.66,4.31)
\psline(6.66,4.31)(7.14,3.7)
\psline(7.14,3.7)(7.05,3.21)
\psline(8.65,5.25)(7.34,5.2)
\psline(4.54,5.16)(4.14,4.49)
\psline(3.6,5.35)(4.54,5.16)
\psline(4.14,4.49)(6.66,4.31)
\psline(1.04,3.21)(1.09,2.3)
\psline(1.09,2.3)(2.14,3.02)
\psline(2.38,3.98)(2.14,3.02)
\psline(2.14,3.02)(2.59,2.77)
\psline(2.59,2.77)(2.8,2.24)
\psline(2.8,2.24)(2.74,1.08)
\psline(2.74,1.08)(3.05,0.28)
\psline(3.05,0.28)(3.83,-1.34)
\psline(3.93,3.27)(2.74,1.08)
\psline(3.96,2.04)(6.64,2.28)
\psline(0.16,1.13)(0.66,-1.36)
\psline(0.66,-1.36)(1.49,-2.25)
\psline(1.49,-2.25)(2.21,-2.53)
\psline(2.21,-2.53)(3.93,-3.55)
\psline(0.04,-0.64)(0.66,-1.36)
\psline(1.16,-3.66)(1.23,-2.84)
\psline(1.23,-2.84)(1.49,-2.25)
\psline(1.04,0.26)(1.69,1.08)
\psline(1.69,1.08)(2.8,2.24)
\psline(2.04,-0.56)(2.72,-1.29)
\psline(2.72,-1.29)(3.05,0.28)
\begin{scriptsize}
\psdots[dotstyle=*,linecolor=blue](6.02,-0.14)
\psdots[dotstyle=*,linecolor=xdxdff](3.33,-5.5)
\psdots[dotstyle=*,linecolor=xdxdff](3.27,-4.9)
\psdots[dotstyle=*,linecolor=xdxdff](3.69,-4.56)
\psdots[dotstyle=*,linecolor=xdxdff](5.08,-4.54)
\psdots[dotstyle=*,linecolor=xdxdff](3.93,-3.55)
\psdots[dotstyle=*,linecolor=xdxdff](5.86,-3.64)
\psdots[dotstyle=*,linecolor=xdxdff](6.62,-3.08)
\psdots[dotstyle=*,linecolor=xdxdff](3.83,-1.34)
\psdots[dotstyle=*,linecolor=xdxdff](4.02,0)
\psdots[dotstyle=*,linecolor=xdxdff](6.32,0.26)
\psdots[dotstyle=*,linecolor=xdxdff](4.62,-5.97)
\psdots[dotstyle=*,linecolor=xdxdff](6.22,-5.64)
\psdots[dotstyle=*,linecolor=xdxdff](7.3,-4.97)
\psdots[dotstyle=*,linecolor=xdxdff](8.02,-4.17)
\psdots[dotstyle=*,linecolor=xdxdff](6.78,-4.07)
\psdots[dotstyle=*,linecolor=xdxdff](7.75,-5.89)
\psdots[dotstyle=*,linecolor=xdxdff](9.94,-4.68)
\psdots[dotstyle=*,linecolor=xdxdff](9.57,-4.34)
\psdots[dotstyle=*,linecolor=xdxdff](9.98,-3.2)
\psdots[dotstyle=*,linecolor=xdxdff](10.33,-1.42)
\psdots[dotstyle=*,linecolor=xdxdff](9.72,-1.65)
\psdots[dotstyle=*,linecolor=xdxdff](9.13,-1.75)
\psdots[dotstyle=*,linecolor=xdxdff](9.02,0)
\psdots[dotstyle=*,linecolor=xdxdff](6.64,2.28)
\psdots[dotstyle=*,linecolor=xdxdff](7.63,-1.33)
\psdots[dotstyle=*,linecolor=xdxdff](11.28,-3.03)
\psdots[dotstyle=*,linecolor=xdxdff](10.96,-2.57)
\psdots[dotstyle=*,linecolor=xdxdff](10.95,-0.99)
\psdots[dotstyle=*,linecolor=xdxdff](12,-0.68)
\psdots[dotstyle=*,linecolor=xdxdff](11.53,2.24)
\psdots[dotstyle=*,linecolor=xdxdff](8.48,4.21)
\psdots[dotstyle=*,linecolor=xdxdff](9.21,3.03)
\psdots[dotstyle=*,linecolor=xdxdff](8.44,3.04)
\psdots[dotstyle=*,linecolor=xdxdff](7.05,3.21)
\psdots[dotstyle=*,linecolor=xdxdff](8.06,2.06)
\psdots[dotstyle=*,linecolor=xdxdff](10.81,1.29)
\psdots[dotstyle=*,linecolor=xdxdff](9.35,2.08)
\psdots[dotstyle=*,linecolor=xdxdff](9.04,1.63)
\psdots[dotstyle=*,linecolor=xdxdff](10.63,3.7)
\psdots[dotstyle=*,linecolor=xdxdff](7.34,5.2)
\psdots[dotstyle=*,linecolor=xdxdff](5.5,4.83)
\psdots[dotstyle=*,linecolor=xdxdff](6.66,4.31)
\psdots[dotstyle=*,linecolor=xdxdff](7.14,3.7)
\psdots[dotstyle=*,linecolor=xdxdff](8.65,5.25)
\psdots[dotstyle=*,linecolor=xdxdff](4.54,5.16)
\psdots[dotstyle=*,linecolor=xdxdff](4.14,4.49)
\psdots[dotstyle=*,linecolor=xdxdff](3.6,5.35)
\psdots[dotstyle=*,linecolor=xdxdff](1.04,3.21)
\psdots[dotstyle=*,linecolor=xdxdff](1.09,2.3)
\psdots[dotstyle=*,linecolor=xdxdff](2.14,3.02)
\psdots[dotstyle=*,linecolor=xdxdff](2.38,3.98)
\psdots[dotstyle=*,linecolor=xdxdff](2.59,2.77)
\psdots[dotstyle=*,linecolor=xdxdff](2.8,2.24)
\psdots[dotstyle=*,linecolor=xdxdff](2.74,1.08)
\psdots[dotstyle=*,linecolor=xdxdff](3.05,0.28)
\psdots[dotstyle=*,linecolor=xdxdff](3.93,3.27)
\psdots[dotstyle=*,linecolor=xdxdff](3.96,2.04)
\psdots[dotstyle=*,linecolor=xdxdff](0.16,1.13)
\psdots[dotstyle=*,linecolor=xdxdff](0.66,-1.36)
\psdots[dotstyle=*,linecolor=xdxdff](1.49,-2.25)
\psdots[dotstyle=*,linecolor=xdxdff](2.21,-2.53)
\psdots[dotstyle=*,linecolor=xdxdff](0.04,-0.64)
\psdots[dotstyle=*,linecolor=xdxdff](1.16,-3.66)
\psdots[dotstyle=*,linecolor=xdxdff](1.23,-2.84)
\psdots[dotstyle=*,linecolor=xdxdff](1.04,0.26)
\psdots[dotstyle=*,linecolor=xdxdff](1.69,1.08)
\psdots[dotstyle=*,linecolor=xdxdff](2.04,-0.56)
\psdots[dotstyle=*,linecolor=xdxdff](2.72,-1.29)
\psdots[dotstyle=*,linecolor=xdxdff](1.23,-2.84)
\end{scriptsize}
\end{pspicture*}
\caption{Tentative description of the discrete radial Web}
\label{fig:figure1}
\end{figure}
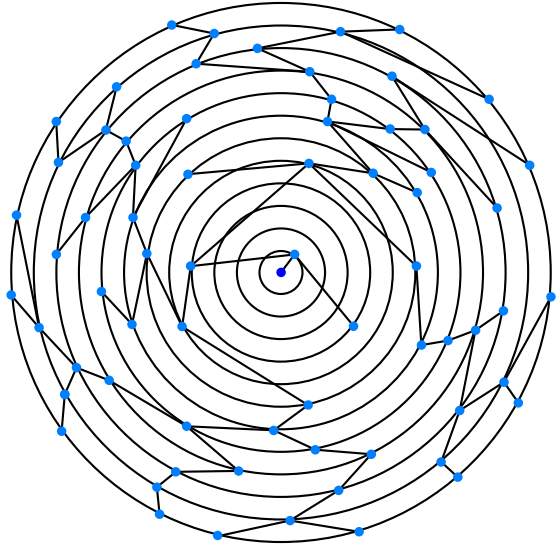

%%%%%%%%%%%%%%%%%%%%%%%%%%%%%%%%%%%%%%%%%%%%%%%%%%%%%%%%%%%%%%%%%%%%%%%%%%%%%%%%%%%%%%%%%%%%%%%%%%%%%%%%%%%%%%%%%%%%%%%%%%%%%%%
%%%%%%%%%%%%%%%%%%%%%%%%%%%%%%%%%%%%%%%%%%%%%%%%%%%%%%%%%%%%%%%%%%%%%%%%%%%%%%%%%%%%%%%%%%%%%%%%%%%%%%%%%%%%%%%%%%%%%%%%%%%%%%%
\section{Definitions, notation and main result} \label{p}
%%%%%%%%%%%%%%%%%%%%%%%%%%%%%%%%%%%%%%%%%%%%%%%%%%%%%%%%%%%%%%%%%%%%%%%%%%%%%%%%%%%%%%%%%%%%%%%%%%%%%%%%%%%%%%%%%%%%%%%%%%%%%%%
%%%%%%%%%%%%%%%%%%%%%%%%%%%%%%%%%%%%%%%%%%%%%%%%%%%%%%%%%%%%%%%%%%%%%%%%%%%%%%%%%%%%%%%%%%%%%%%%%%%%%%%%%%%%%%%%%%%%%%%%%%%%%%%
Let $z=(x,y) \in \mathbb{R}^2$ with $y<0$. We will use the letter $\vartheta_z$ to denote the angle between $z$ and $(0,-1)$. Let 
$\tilde{\vartheta}_z=sgn(x)\vartheta_z$. Let $\alpha\in(0,1)$ and $\delta\in(1/4,1/3)$ be fixed 
(for technical reasons). For each $n\in\bbN$, define $\theta_n=n^{\delta/2-1/2}$, $\varphi_n=n^{\delta-1/2}$ and

\begin{equation}
 \mathcal{B}_n:=\mathcal{B}_n(\alpha,\delta)=\left\{ z\in \bbR\times [-n,-n\alpha] : |\tilde{\vartheta}_z|
 \leq \theta_n/2, |z|\leq n  \right\}, \nonumber
\end{equation}
and
\begin{equation}
 \mathcal{A}_n:=\mathcal{A}_n(\alpha,\delta)=\left\{ z\in \bbR\times [-n,-n\alpha] : |\tilde{\vartheta}_z|
 \leq \varphi_n/2, |z|\leq n  \right\}. \nonumber
\end{equation}
Let $\{ P_n \}_{n=1}^\infty$ be a sequence of independent Poisson processes with rate $1$ on $\mathbb{R}$. Then
\begin{equation}\label{process}
\hat{P}_n:=\{ s\in\mathbb{R}^2 : s = n\left(\sin\left( x/n \right), -\cos (x/n)  \right)\hspace{.2cm}
\text{for some $x\in P_n\cap [-n\pi,n\pi]$} \}
\end{equation}
is a Poisson process with rate $1$ on $\{ s\in\mathbb{R}^2: \|s \|=n \}$ where $\|.\|$ stands for the euclidean norm. For each $n\in \mathbb{N}$, let $\varGamma_n:=\hat{P}_n\cap \mathcal{B}_n$. Now fix some $\kappa \in (0,1/2-\delta)$ and $k \in \{n,n-1,...,\lfloor n\alpha  \rfloor-1 \}$. Now we define a random path starting at $z \in \varGamma_k$ as follows:

\begin{enumerate}
 \item 
 If $\hat{P}_{k-1}\cap \mathcal{A}_n=\emptyset$, then join $z$ and the origin with a straight line.
 \item 
 If $\hat{P}_{k-1}\cap \mathcal{A}_n\neq\emptyset$, let $z_1$ be the point in $\hat{P}_{k-1}\cap \mathcal{A}_n$ minimizing the distance to $z$. 
 \begin{enumerate}
 \item 
 If $z=(x,y), z_1=(x_1,y_1)$ and $|x-x_1|>n^{\kappa}$, then join $z$ and the origin with a straight line.
 
 \item 
If $|x-x_1| \leq n^{\kappa}$, then join $z$ and $z_1$ with a straight line.

\end{enumerate}
\end{enumerate}
If the origin has not been reached, then repeat this procedure until the origin is reached. Denote by  $\gamma_z$ the random trajectory 
obtained in this way. In order to figure out this definition see figure~{\ref{fig:representation}}.

Let 

\begin{equation}
\gamma^n:=\left\{ \gamma_z : z\in \Gamma_k\hspace{.2cm}\text{for some} k\in\{n,n-1,...,
\lfloor n\alpha\rfloor\} \right\}
\end{equation}
be the set of all coalescence paths starting in $\mathcal{B}_n$. The random set of paths $\left(\gamma^n\right)_n$ is called 
the {\it{Discrete Radial Poissonian Web}}. See figure~{\ref{fig:representation}} for a schematic representation of two coalescing random paths.

\psset{unit=1.3}
\definecolor{xdxdff}{rgb}{0.0, 0.5, 1.0}.
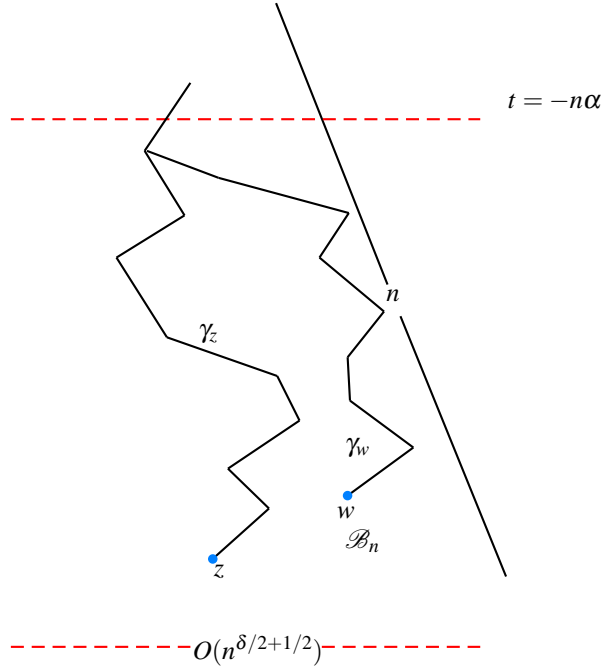
\begin{figure}[h]
\begin{pspicture*}(-2,-5.5)(15,6)
%\pscustom[linecolor=blue,fillcolor=blue,fillstyle=solid,opacity=0.1]{
%\parametricplot[plotstyle=curve, linewidth=2pt]{
%4.331882603272325}{5.092895357497055}{1*10.77*cos(t)+0*10.77*sin(t)+6|0*10.77*cos(t)+1*10.77*sin(t)+6}
%\lineto(6,6)\closepath}
\rput[tl](5.1,-4.8){$O(n^{\delta/2+1/2})$}
\pswedge[fillstyle=solid,fillcolor=xdxdff!10](6,5.95){10.5}{248}{292}
\psline(6.52,5.98)(8.42,1.18)
\psline(8.64,0.64)(10.44,-3.8)
\rput[tl](5.8,5.6){}
\rput[tl](8.40,1.1){$ n $}
\rput[tl](7.7,-3.0){$ \mathcal{B}_n$}
\psline[linestyle=dashed](2,4)(10,4)
\psline[linestyle=dashed](2,-5)(4.9,-5)
\psline[linestyle=dashed](8,-5)(10,-5)
\rput[tl](10.45,4.08){$ t= -  n\alpha $}
\rput[tl](5.46,-3.6){$ z $}
\psline(5.44,-3.5)(6.4,-2.64)
\psline(6.4,-2.64)(5.7,-1.96)
\psline(5.7,-1.96)(6.92,-1.14)
\psline(6.92,-1.14)(6.54,-0.38)
\psline(6.54,-0.38)(4.66,0.28)
\psline(4.66,0.28)(3.8,1.64)
\psline(3.8,1.64)(4.96,2.36)
\psline(4.96,2.36)(4.28,3.46)
\psline(4.28,3.46)(5.06,4.62)
\rput[tl](5.2,0.5){$ \gamma_z $}
\psline(7.74,-2.42)(8.86,-1.6)
\psline(8.86,-1.6)(7.78,-0.8)
\psline(7.78,-0.8)(7.74,-0.06)
\psline(7.74,-0.06)(8.36,0.72)
\psline(8.36,0.72)(7.26,1.64)
\psline(7.26,1.64)(7.76,2.4)
\psline(7.76,2.4)(5.54,3)
\psline(5.54,3)(4.32,3.46)
\rput[tl](7.56,-2.6){$ w $}
\rput[tl](7.7,-1.44){$ \gamma_w $}
\begin{scriptsize}
\psdots[dotstyle=*,linecolor=xdxdff](5.44,-3.5)
\psdots[dotstyle=*,linecolor=xdxdff](7.74,-2.42)
% \psdots[dotstyle=*,linecolor=red](7.74,-2.42)
\end{scriptsize}
\end{pspicture*}
\caption{Schematic representation of two random paths in the DRPW.}
\label{fig:representation}
\end{figure}
 
\begin{rem}
\label{1}
We explain briefly why we introduce conditions $(1)$ and $(2)(a)$ above. 
%First, note that these conditions guarantee that the random paths thus obtained are functions. In order to do that, 
We require $\delta \in (0,1/3)$ in order to guarantee that, with high probability, the family of random 
polygonal curves starting from $\mathcal{B}_n$ are functions. Also, $\delta > 1/4$ and 
$\alpha \in (0,1)$ implies that, with high probability, this family of curves is contained in 
$\mbox{int} (\mathcal{A}_n)$. Then note that the events defining conditions $(1)$ and $(2)(a)$ occur with probability going to zero exponentially fast. These conditions impose modifications in the random paths described in the tentavive description given in the introduction. The properties just described say that these modifications are not significant when $n$ is large enough.
\end{rem}

\subsection{The Brownian Web: characterization}
\label{TBW:C}
First, we give a brief description of the BW which follows closely the description given 
in \cite{ss}, see also \cite{FINR} and the appendix in \cite{ss1}. Consider the extended 
plane $\bar{\mathbb{R}}^2 = [-\infty,\infty]^2$ as the completion of $\mathbb{R}^2$ under the metric
\[
\rho((x_1,t_1),(x_2,t_2)) = |\tanh(t_1) - \tanh(t_2) | \vee \Big| \frac{\tanh(x_1)}{1+|t_1|} - \frac{\tanh(x_2)}{1+|t_2|} \Big| \, ,
\]
and let $\bar{\rho}$ be the induced metric on $\bar{\mathbb{R}}^2$. 
In $(\bar{\mathbb{R}}^2,\bar{\rho})$, the lines $[-\infty,\infty] \times \{ \infty \}$ 
and $[-\infty,\infty] \times \{ - \infty \}$ correspond respectively to single points 
$(\star, \infty)$ and $(\star,-\infty)$, see picture 2 in \cite{ss1}. Denote by $\Pi$ the set of all continuous 
paths in $(\bar{\mathbb{R}}^2,\bar{\rho})$ of the form 
$\pi: t \in [\sigma_\pi,\infty] \rightarrow (f_\pi(t),t) \in (\bar{\mathbb{R}}^2,\bar{\rho})$ for some 
$\sigma_\pi \in [-\infty,\infty]$ and $f_\pi : [\sigma_\pi,\infty] \rightarrow [-\infty,\infty] \cup \{\star \}$. 
For $\pi_1$, $\pi_2 \in \Pi$, define $d(\pi_1,\pi_2)$ by
\[
|\tanh(\sigma_{\pi_1}) - \tanh(\sigma_{\pi_2}) | \vee \sup_{t \ge \sigma_{\pi_1} \wedge \sigma_{\pi_2}}
\Big| \frac{\tanh(f_{\pi_1}(t \vee \sigma_{\pi_1}))}{1+|t|} - \frac{\tanh(f_{\pi_2}(t \vee \sigma_{\pi_2}))}{1+|t|} \Big| \, .
\]
Thus, we have a metric in $\Pi$ such that $(\Pi,d)$ is a complete separable metric space. 
Now define $\mathcal{H}$ as the space of compact sets of $(\Pi,d)$ with the topology 
induced by the Hausdorff metric. Then $\mathcal{H}$ is a complete separable metric space. 
The Brownian web is a random element $\mathcal{W}$ of $\mathcal{H}$ whose distribution is uniquely 
characterized by the following three properties (see Theorem 2.1 in \cite{FINR}):
\begin{enumerate}
\item[(a)] For any deterministic $z \in \mathbb{R}^2$, almost surely there is a unique path 
$\pi_z$ of $\mathcal{W}$ that starts at $z$.
\item[(b)] For any finite deterministic set of points $z_1$, ... ,$z_k$ in $\mathbb{R}^2$, 
the collection $(\pi_1,...,\pi_n)$ is distributed as coalescing Brownian motions independent 
up to the time of coalescence.
\item[(c)] For any deterministic countable dense subset $\mathcal{D} \subset \mathbb{R}^2$, 
almost surely, $\mathcal{W}$ is the closure of $\{ \pi_z : z \in \mathcal{D} \}$ in $(\Pi,d)$.
\end{enumerate}

\subsection{The restricted Brownian Web}\label{TBRR}

Let $r,s\in\mathbb{R}$ be such that $r<s$ and let $A_{r,s}=\mathbb{R}\times[r,s]$. 
Let $\bar{A}_{r,s}$ be the compactification of $A_{r,s}$ under the metric
$\rho$.
\noindent $\bar{A}_{r,s}$ may be thought as the image of $[-\infty,\infty]\times[r,s]$ under the mapping

\begin{equation}
(x,t)\rightsquigarrow \left(\Phi (x,t),\Psi (t)\right) \equiv \left(\frac{\tanh (x)}{1 +
|t|},\tanh (t) \right). \nonumber
\end{equation}
For $t_0\in[r,s]$, let $\mathcal{C}_{r,s}[t_0]$ be the set of continuous functions $f$ 
from $[t_0,s]$ to $[-\infty,\infty]$. Define

\begin{equation}
\Pi_{r,s}=\bigcup_{t_0\in[r,s]}\mathcal{C}_{r,s}[t_0]\times\{ t_0  \} \nonumber
\end{equation}

\noindent where $(f,t_0)\in\Pi_{r,s}$ represent a path in $\bar{A}_{r,s}$ starting at $(f(t_0),t_0)$.
Denote by $\hat{f}$ the function that extends $f$ to all $[r,s]$ by setting it equal 
to $f(t)=f(t_0)$ for $r\leq t\leq t_0$. Then, define

\begin{equation}
d_{r,s}((f_1,t_1),(f_2,t_2))=\sup_{r\leq t\leq s}|\Phi(\hat{f}_1(t),t) - \Phi(\hat{f}_2(t),t)| \vee
|\Psi (t_1) - \Psi (t_2)|. \nonumber
\end{equation}

Thus, $(\Pi_{r,s},d_{r,s})$ is a complete, separable metric space.

\paragraphmark

Let $\mathcal{H}_{r,s}$ be the set of compact subsets of 
$(\Pi_{r,s},d_{r,s})$ with $d_{\mathcal{H}_{r,s}}$ the induced Hausdorff metric given by

\begin{equation}\label{aches}
d_{\mathcal{H}_{r,s}}(K_1,K_2)=\sup_{g_1\in K_1}\inf_{g_2\in K_2}d_{r,s}(g_1,g_2)\vee 
\sup_{g_2\in K_2}\inf_{g_1\in K_1}d_{r,s}(g_1,g_2).
\end{equation}
Then, $(\mathcal{H}_{r,s}, d_{\mathcal{H}_{r,s}})$ is also a complete, separable metric space. 
Also, denote by $\mathcal{F}_{\mathcal{H}_{r,s}}$ the $\sigma$-algebra associated to $d_{\mathcal{H}_{r,s}}$.

\paragraphmark

Now, we define and characterize the restricted Brownian Web in terms of the standard Brownian Web. Also, we give a 
criteria of convergence to the restricted web.

Let $r,s\in\mathbb{R}$ be fixed such that $r<s$. Let

\[
\hat{\mathcal{H}}_r^s:=\left\{ \mathcal{J}\in \mathcal{H} : 
\text{There is} \hspace{.2cm} (f,t_0)\in \mathcal{J} \hspace{.2cm}\text{with}\hspace{.2cm} t_0\in [r,s]   \right\}.
\]
In other words, $\hat{\mathcal{H}}_r^s$ is a subset of $\mathcal{H}$ formed by compact families of 
trajectories having at least one trajectory starting from 
$\bar{A}_{r,s}$. Since $\hat{\mathcal{H}}_r^s$ is a closed subset of $\mathcal{H}$ it 
follows that $\hat{\mathcal{H}}_r^s$ is also a complete, separable metric space.

Consider the mapping  $\mathfrak{T}$ from $\hat{\mathcal{H}}_r^s$ to $\mathcal{H}_{r,s}$ defined by

\begin{equation}
\mathfrak{T}(K)=\{ (f,t_0)\in K \hspace{.1cm}\text{restricted to the set $\bar{A}_{r,s}$}: t_0\in [r,s]\}. \nonumber
\end{equation}
In other terms, $\mathfrak{T}(K)$ is the set of trajectories of $K$ starting at points in $\bar{A}_{r,s}$ 
and restricted to the set $\bar{A}_{r,s}$. We claim that $\mathfrak{T}$ is well defined. Indeed, 
for $K\in \hat{\mathcal{H}}_r^s$ set

$$
K^{\prime} = \left\{(f,t_0)\in K : t_0 \in [r,s] \right\}.
$$

\noindent Then, $K^{\prime}$ is a closed subset of $K$. Therefore, it is a compact subset of trajectories. 
Now, it follows from this that the set of trajectories of $K^{\prime}$ restricted to the set $\bar{A}_{r,s}$ is compact, i.e
$\mathfrak{T}(K)=\mathfrak{T}(K^{\prime})$ is compact. This proves the claim.

Observe that $K_n\underset{n\rightarrow\infty}{\longrightarrow}K$ implies 
$\mathfrak{T}(K_n)\underset{n\rightarrow\infty}{\longrightarrow}\mathfrak{T}(K)$. Then, we may conclude that 
the mapping $\mathfrak{T}$ is continuous. Next, we characterize the {\it{restricted Brownian Web}}.

\begin{thm} \label{restricta1} 
There exists a $(\mathcal{H}_{r,s}, \mathfrak{F}_{\mathcal{H}_{r,s}})$-valued random variable 
$\bar{\mathcal{W}}_{r,s}$ whose distribution is uniquely determined by the following properties.

\begin{enumerate}
	\item[(a) ] From any deterministic point $(x,t)\in A_{r,s}$, there is almost surely a unique path $W_{x,t}$ starting from $(x,t)$.
	\item[(b)] For any deterministic $n, (x_1,t_1),...,(x_n,t_n)\in A_{r,s}$ the joint distribution 
	of $W_{x_1,t_1},...,W_{x_n,t_n}$ is that of coalescing Brownian motions with unit diffusion constant.
	\item[(c)] For any deterministic, dense countable subset $\mathcal{D}$ of $A_{r,s}$, 
	almost surely, $\bar{{\mathcal{W}}}_{r,s}$ is the closure in \newline $(\mathcal{H}_{r,s}, 
	\mathfrak{F}_{\mathcal{H}_{r,s}})$ of $\{W_{x,t}:(x,t)\in\mathcal{D}   \}$.
\end{enumerate}
\end{thm}

\begin{proof}
From Theorem 2.1 in \cite{FINR} we know that, almost surely, $\bar{\mathcal{W}}_{r,s}\in\hat{\mathcal{H}}_r^s$. 
Then, the result follows from the continuity of $\mathfrak{T}$ and Theorem 2.1 in \cite{FINR} 
(See properties $(a), (b)$ and $(c)$ in subsection \ref{TBW:C}).
\end{proof}

\begin{defn}
We call the random variable $\bar{\mathcal{W}}_{r,s}$, the restricted Brownian web.
\end{defn}

\subsection{Restricted Brownian Web: convergence criteria}
Now we give the convergence criteria to the restricted Brownian Web. Let $\mathcal{D}$ be a 
countable dense subset of $\mathbb{R}^2$. Also, let $t_0\in[r,s]$, $t>0$, $a<b$ and 
$\mathcal{M}$ be a $(\mathcal{H}_{r,s},\mathfrak{F}_{\mathcal{H}_{r,s}})$-valued random variable. 
Let $\eta_{\mathcal{M}}(t_0,t;a,b)$ be the $\{0,1,2,...\}$-valued random variable representing 
the number of distinct points in $\mathbb{R}\times\{t_0+t \}$
which are touched for some paths from $\mathcal{M}$ which also cross 
the segment $[a,b]\times\{ t_0 \}$. If $t_0+t>s$, make $\eta_{\mathcal{M}}(t_0,t;a,b)\equiv 0$.

Now we state the convergence criteria to the restricted Brownian web. Let $\mathcal{D}$ be countable dense subset of $A_{r,s}$.

\begin{thm}\label{omega4}
Suppose that $\mathcal{X}_1,\mathcal{X}_2,...$ are $(\mathcal{H}_{r,s},\mathfrak{F}_{\mathcal{H}_{r,s}})$-valued 
random variables with non crossing paths. If the following three conditions are valid, the distribution 
$\mu_n$ of $\mathcal{X}_n$ converges to the distribution $\mu_{\bar{\mathcal{W}}_{r,s}}$ of the restricted Brownian web.

\begin{itemize}
\item $(I)$ For any deterministic $y_1,y_2,...,y_m \in \mathcal{D}$, there exist $\theta_n^{y_1},...,
\theta_n^{y_m} \in \mathcal{X}_n$ such that $\theta_n^{y_1},...,\theta_n^{y_m}$ converge in distribution as 
$n \rightarrow \infty$ to coalescing Brownian motions (with unit diffusion constant) starting at $y_1,y_2,...,y_m${\color{red}{,}}
\item $(B_1)$ $\forall \ t > 0, \limsup_{n \rightarrow \infty} \sup _{(a,t_0) \in A_{r,s}} \mu_n (
{\eta}_{\mathcal{X}_n}(t_0,t;a,a+\epsilon) \geq 2) \rightarrow 0$ as $\epsilon \rightarrow 0+${\color{red}{,}}
\item $(B_2)$ $\forall \ t > 0, \epsilon^{-1} \limsup_{n \rightarrow \infty} \sup _{(a,t_0) \in A_{r,s}} \mu_n (
{\eta}_{\mathcal{X}_n}(t_0,t;a,a+\epsilon) \geq 3) \rightarrow 0$ as $\epsilon \rightarrow 0+${\color{red}{.}}
\end{itemize}
\end{thm}

\noindent {\bf Sketch of the Proof} The fact that the paths do not cross, along with condition $(I)$ 
imply tightness of the sequence $(\mathcal{X}_n)_{n=1}^{\infty}$. Conditions $(B_1)$ and $(I)$ guarantee 
that any subsequence limit $\mathcal{X}$ of $\{ \mathcal{X}_n  \}$ contains a version of the restricted Brownian web. 
It remains to argue that$\mathcal{X}$ does not contain anything else. Conditions $(B_1)$ and $(B_2)$ together imply that 
$\mathbb{E}[\eta_{\mathcal{X}}(t_0,t;a,b)] \leq \mathbb{E}[\eta_{\overline{\mathcal{W}}}(t_0,t;a,b)] = 
1+\frac{b-a}{\sqrt{\pi t}}$ for any $t_0 \in [r,s], t>0, a,b \in \mathbb{R}$. Then, Theorem $4.6$ in 
\cite{FINR} implies that $\mathcal{X}$ contains nothing else than paths from $\bar{\mathcal{W}}_{r,s}$.

\subsection{The T-Brownian Web and statement of the main result}
The main result of this work is to show that under the diffusive scaling,
$\gamma^n$ converges in distribution to a continuous mapping of the restriction of the 
Brownian web to the set $\mathbb{R}\times[0,1/\alpha-1]$, where $\alpha \in (0,1)$ is fixed.

Unless otherwise stated, from now on  $\alpha$ will be any fixed number in $\left(0,1\right)$. Let $F_{\alpha}:=\mathbb{R}\times[-1,-\alpha]$ and $G_\alpha:=\mathbb{R}\times[0,1/\alpha-1]$. Let
$\psi : F_\alpha \mapsto G_\alpha$ be the homeomorphism defined by

\[
\psi(x,t)=\left( \frac{x}{|t|}, \frac{1}{|t|}-1  \right).
\]
Then,

\[
\psi^{-1}(x^{\prime},t^{\prime})=\left( \frac{x^{\prime}}{t^{\prime}+1},-\frac{1}{t^{\prime}+1}   \right).
\]
The map $\psi$ is an homeomorphism between $\mathbb{R} \times [-1,-\alpha]$ and $\mathbb{R} \times [0,1/\alpha - 1]$. Loosely speaking, $\psi$ maps each trajectory of the Brownian Bridge Web into a trajectory of the Brownian Web.

Let $r,s\in\mathbb{R}$ be such that $r<s$. Denote by $\mathcal{H}_{r,s}$ the set 
of compact subsets of trajectories starting from points inside $\mathbb{R}\times[r,s]$ and 
restricted to this same set under the metric defined in (\ref{aches}).

Let $T: \mathcal{H}_{-1,-\alpha}\longrightarrow\mathcal{H}_{0,1/\alpha-1}$ be the map defined by

\begin{equation}
T(\mathcal{G})=\left\{(\psi(f(t),t)) \right\}_{(f(t),t)\in \mathcal{G}}.
\end{equation}
%\label{transformacion}
%Then,
%
%\[
%T^{-1}(\mathcal{F})=\left\{(\psi^{-1}(f(t),t)) \right\}_{(f(t),t)\in \mathcal{F}}.
%\]
For $\mathcal{F}\in \mathcal{H}_{-n,-n\alpha}$ define

\begin{equation}
G_n[\mathcal{F}]=\left\{ \left(\frac{f(t)}{\sqrt{n}},\frac{t}{n} \right): (f(t),t)\in \mathcal{F}  \right\}. \nonumber
\end{equation}
and note that $G_n[\mathcal{F}]\in \mathcal{H}_{-1,-\alpha} $. 

By abuse of notation denote by $\hat{\gamma}^n$ the set of random paths $\left(\gamma^n\right)_n$ restricted to the set $\mathbb{R}\times [-n,-n\alpha]$. Then let  $\Phi^n:=G_n[\hat{\gamma}^n]$ be the restricted DRPW under the diffusive scaling and note that $\Phi^n\in\mathcal{H}_{-1,-\alpha}$. Our main result is the following:

\begin{thm} \label{omega}
Let $\alpha \in \left(0,1\right)$ be fixed. Then the rescaled discrete radial Poissonian web $\Phi^n$ converges in distribution to $T^{-1}(\bar{\mathcal{W}}_{0,1/\alpha-1})$ as $n \rightarrow \infty$.
\end{thm}

\begin{rem} 
We observe that to prove the last statement it suffices to show that $T(\Phi^n)$ converges in distribution to $\bar{\mathcal{W}}_{0,1/\alpha-1}$.
\end{rem}

\section{Proof of the main result} \label{sec:auxresults}
For convenience, along this section we will remind the reader which paths are macroscopic and which paths are microscopic. By a macroscopic path we mean a trajectory in the Discrete Radial Web and by a microscopic path we mean a path in the Discrete Radial Web under diffusive scaling.
\subsection{Weak convergence of $\Phi^n$}

In order to prove Theorem \ref{omega} we introduce a family of coalescing random paths approximating the family $\Phi^n$. Since $\Phi^n$ lies inside the interior of the set $\mathcal{A}_n$, the random points used in the linear interpolation used to define $\gamma^n$ may be represented in polar coordinates $\left(r\sin(\theta), -r \cos(\theta)\right)$, where $r\leq n$ and $|\theta|\leq \varphi_n/2$. Therefore, for $n$ large enough we have

\begin{equation}\label{approx}
\left( \frac{r\sin(\theta)}{\sqrt{n}}, \frac{-r \cos(\theta)}{n}  \right)\sim \left( \frac{r\theta}{\sqrt{n}}, -\frac{r}{n}  \right),
\end{equation}
where $\left(a^1_n,a^2_n\right) \sim \left(b^1_n,b^2_n\right)$ means that $a^i_n/b^i_n \rightarrow 1$ as $n \rightarrow \infty$ for $i=1,2$.

Now we build a family of coalescing random paths approximating $\Phi^n$. For simplicity assume that, for $n$ large enough, there exists a path in $\Phi^n$ starting at $(0,-n)$. Then choose the nearest point to $\left(0,-n\right)$ belonging to $\hat{P}_{n-1}\cap \mathcal{A}_n$, which exists with high probability. It follows from (\ref{process}) that this point may be represented as follows:

\begin{equation}
\zeta_1:=(n-1)\left(\sin\left( \frac{x_{n-1}}{n-1}\right), - \cos\left( \frac{x_{n-1}}{n-1}\right)\right), \nonumber
\end{equation}
where $x_{n-1}\in P_{n-1}$. Assuming that we have succeeded in finding the point $\zeta_1$ we look for the point in
%Buscamos ahora el punto que minimize la distancia al punto $\zeta_1$ en 
$\hat{P}_{n-2}\cap \mathcal{A}_n$ minimizing the distance to $\zeta_1$ which exists with high probability. It follows from (\ref{process}) that this point may be represented as follows:

\begin{equation}
\zeta_2 := (n-2)\left(\sin\left( \frac{x_{n-2}}{n-2}\right), - \cos\left( \frac{x_{n-2}}{n-2}\right)\right), \nonumber
\end{equation}
where $\frac{x_{n-2}}{n-2}\in \frac{1}{n-2}P_{n-2}$. Continuing in this way, if possible, we get a point $\zeta_k$ which may be represented as follows:

\begin{equation}
\zeta_k:= (n-k)\left(\sin\left( \frac{x_{n-k}}{n-k}\right), - \cos\left( \frac{x_{n-k}}{n-k}\right)\right). \nonumber
\end{equation}
Now, for each $k \in \mathbb{N}$, let $\omega_k:=\arg\min\{ |x|: x\in P_k \}$. Then, given the events $(1)$, $(2)$ and using the traslation invariance of a Poisson process, it is not difficult to show that $\zeta_k$ is equally distributed with the random point 

\begin{equation}
(n-k)\left( \sin\left( \sum_{j=1}^k \frac{\omega_j}{n-j}\right),  -\cos\left(  \sum_{j=1}^k \frac{\omega_j}{n-j}\right)\right). \nonumber
 \end{equation}

It follows from (\ref{approx}) that

\begin{equation}  \label{refcorr}
(n-k)\left( \frac{\sin\left( \sum_{j=1}^k \frac{\omega_j}{n-j}\right)}{\sqrt{n}}, -\frac{\cos\left(  \sum_{j=1}^k \frac{\omega_j}{n-j}\right)}{n}\right)\sim
\left( \frac{n-k}{n}\right)\left( \sqrt{n}\sum_{j=1}^k \frac{\omega_j}{n-j}, -1 \right)=\xi_k. 
\end{equation}

Let $z \in\mathbb{R}\times(-\infty,0)$ and let $r \left(\sin\theta,-\cos\theta \right)$ be its representation in polar coordinates for some $r >0$ and $\theta \in \left(-\pi/2,\pi/2\right)$. Define a function $\Lambda :\mathbb{R}\times(-\infty,0)\longrightarrow \mathbb{R}^2$ by the formula $\Lambda(z)=(r\theta,-r)$. Denote this new family of random coalescing paths by

\begin{equation*}
\xi^n:=\Lambda[\Phi^n]:=\left\{ \Lambda(z) : z\in \Phi^n\right\}
\end{equation*}
Denote by $\xi_z^n$ a path in $\xi^n$ starting from $z$ and note that $\xi^n$ is a family of random macroscopic paths.

\begin{prop}\label{prop1}
Let $\rho$ denote the Haussdorff metric. Then $\rho\left( \xi^n, \Phi^n \right)\to 0$ as $n\to\infty$, in probability.
\end{prop}
It follows from Proposition \ref{prop1} that in order to prove Theorem \ref{omega} it suffices to show that

\begin{equation*}
\xi^n\overset{\mathcal{D}}{\underset{n\to\infty}{\longrightarrow}} T^{-1}\left(\bar{\mathcal{W}}_{0,1/\alpha -1}\right), 
\end{equation*}
or the equivalent statement

\begin{equation*}
T\left( \xi^n \right)\overset{\mathcal{D}}{\underset{n\to\infty}{\longrightarrow}} \bar{\mathcal{W}}_{0,1/\alpha -1}
\end{equation*}
Therefore, we must study the behavior, in distribution, of the random set of microscopic coalescing paths $T(\xi^n)$. 

\begin{lem}\label{1path}
Let $\Delta_n := \left\{ \left( \frac{x}{\sqrt{n}}, \frac{y}{n} \right): (x,y)\in \cup_{k=1}^n \Gamma_n \right\}$ and let $(x_0,t_0)\in\mathbb{R}\times[-1,-\alpha]$ be fixed. Then, almost surely, $\Delta_n$ is a compact set and
\begin{equation}
\mathbb{P}\left[  \rho\left( \Delta_n , \{ (x_0,y_0) \}\right)> 2n^{\delta /2-1/2}  \hspace{.1cm} i.o\right]=0. \nonumber
\end{equation}
\end{lem}

\begin{proof}
It follows from the properties of Poisson processes and the Borel-Cantelli lemma.
\end{proof}
The previous lemma tell us that given a deterministic point $(x_0,y_0)\in\mathbb{R}\times[-1,-\alpha]$, then with high probability  and if $n$ is large enough there exists a trajectory in $\xi^n$ such that its starting point is as close as we desire from $(x_0,y_0)$.

\begin{lem}\label{path}
Let $(x_0,t_0)\in\mathbb{R}\times[-1,-\alpha]$ be fixed and let $\{W_t\}_{t \geq 0}$ be a standard Brownian motion starting at the origin. Then, there exists a sequence $\xi_z^n\in\xi^n$ such that

\begin{equation}
\left\{\xi_z^n\right\}\overset{\mathcal{D}}{\underset{n\rightarrow\infty}{\longrightarrow}}\left\{ t\left(\frac{x_0}{t_0}+\sqrt{2}W_{g(t)}  \right)\right\}_{t\in [t_0, -\alpha]},  \nonumber
\end{equation}
where $W_{g(t)}\sim N(0,g(t))$, $t$ is fixed and $g(t)=\int_{t_0}^{t}\frac{1}{x^2}dx=\frac{1}{|t|}-\frac{1}{|t_0|}$. Also, \hfill \break  $\left\{W_{g(t)}\right\}_{t\in[t_0,-\alpha]}$ has independent increments.
\end{lem}

For the proof of Lemma \ref{path} see the first part of Appendix \ref{sec:App}. Also, and under the same conditions of Lemma \ref{path} above, we have that

\[
T(\xi_z^n)\overset{\mathcal{D}}{\underset{n\rightarrow\infty}{\longrightarrow}}\left\{  x_0^{\prime}+\sqrt{2}W_{t^{\prime}-t_0^{\prime}}  \right\}_{t_0^{\prime}\leq t^{\prime}\leq \frac{1}{\alpha}-1}
\]
where $x_0^{\prime}=\frac{x_0}{|t_0|}$ and $t_0^{\prime}=\frac{1}{|t_0|-1}$.

\subsection{Alternative model}

In order to establish weak convergence of $T\left(\xi^n\right)$ we introduce an alternative model which is tractable and whose weak limit coincides with that of $T\left(\xi^n\right)$. Indeed, this claim easily follows from (\ref{refcorr}) and the definition of $T$. 

For $n$ large enough and $j=0,1,...,\left\lfloor n(1-\alpha)\right\rfloor+1$, let $Q_j=\frac{\sqrt{n}}{n-j}P_j$. If $k_n=\left\lfloor n(1-\alpha)\right\rfloor$, then
$\{Q_j\}_{j=0}^{k_n+1}$ is a family of Poisson process where $Q_j$ has rate $\frac{n-j}{\sqrt{n}}$. The choise of the parameters for this family of point processes deserves some words of explanation. It follows from \ref{refcorr} that trajectories from $\xi^n$ behaves asymptotically, in distribution, as the trajectories obtained after linearly interpolating the points
\begin{equation} \label{ab}
\frac{n-k}{n} \left(\sqrt{n} \displaystyle \sum_{j=1}^k \frac{\omega_j}{n-j},-1\right)
\end{equation}
where $k=1,2, \ldots, \lfloor n\left(1-\alpha\right)\rfloor + 1$. Now observe that $T$ maps any trajectory passing through points of the form \ref{ab} into a trajectory passing through points of the form below
\[
\left(\sqrt{n} \displaystyle \sum_{j=1}^k \frac{\omega_j}{n-j},\frac{k}{n-k}\right).
\]
Finally note that the first coordinate of the vector above behaves as a random walk. Thus the parameters are chosen in order to understand the asymptotic behavior of this random walk.

Now consider the set of random points $S^{(n)}=\cup_{j=0}^{k_n+1}Q_j\times \left\{  \frac{j}{n-j}\right\}$. Then, join with a straight line each point 
$u$ in $Q_j\times \left\{ \frac{j}{n-j} \right\}$ to the nearest point $Q_{j+1}\times\left\{ \frac{j+1}{n-j-1} \right\}$ which exists with probability one. Then repeat this procedure up to time $k_n+1$. Denote by $\beta_u^n$ the trajectory thus obtained. Let $\beta^n:=\left\{\beta_u^n: u\in S^{(n)} \right\}$ denotes the corresponding family of macroscopic coalescing random paths. In order to figure out the behavior of this family of random paths see figure~\ref{fig:representationAM}. Also, let $\mu_n$ the distribution induced by $\beta^n$

\psset{unit=0.7}
\begin{figure}[h]
\begin{pspicture*}(-6.0,-10)(17,8)
\psline(-2,-8)(13,-8)
\psline(-2.01,-7.2)(12.99,-7.2)
\psline(-2,-6)(13,-6)
\psline(-1.95,-3)(13.05,-3)
\psline(-1.98,-1.38)(13.02,-1.38)
\psline(-2,1.03)(13,1.03)
\psline(-1.97,2.97)(13.03,2.97)
\psline(-0.77,-7.24)(0,-6)
\psline(0.4,-8)(0.33,-7.2)
\psline(0.33,-7.2)(0,-6)
\psline(1.2,-7.2)(0,-6)
\psline(3.48,-8)(3.99,-7.2)
\psline(3.99,-7.2)(3.7,-6)
\psline(4.3,-8)(3.99,-7.2)
\psline(5.44,-8)(3.99,-7.2)
\psline(8.31,-8)(7.5,-7.2)
\psline(7.5,-7.2)(8.18,-5.93)
\psline(9.39,-7.2)(9,-6)
\psline(12,-8)(9.39,-7.2)
\rput[tl](13.08,-7.48){$t=0$}
\rput[tl](13.08,-6.34){$t=\frac{1}{n-1}$}
\rput[tl](13.06,-5.18){$t=\frac{2}{n-2}$}
\psline(-1.42,-3)(-0.14,-1.38)
\psline(-0.47,-3)(-0.14,-1.38)
\psline(0.42,-3)(-0.14,-1.38)
\psline(4.73,-3)(6.02,-1.38)
\psline(5.48,-3)(6.02,-1.38)
\psline(7.65,-3)(6.99,-1.38)
\psline(8.25,-3)(6.99,-1.38)
\psline(12.24,-3)(10.29,-1.38)
\rput[tl](13.19,-2.16){$t=\frac{k}{n-k}$}
\rput[tl](13.08,-0.54){$t=\frac{k+1}{[n-(k+1)]}$}
\rput[tl](13.06,1.84){$t=\frac{k_n}{n-k_n}$}
\rput[tl](-4.19,2.55){$t=\frac{1-\alpha}{\alpha}$}
\psline(-2.03,1.74)(12.97,1.78)
\psline(-1.27,1.03)(-0.5,2.97)
\psline(0,1.03)(-0.5,2.97)
\psline(0.64,1.03)(-0.5,2.97)
\psline(3.16,1.03)(2.41,2.97)
\rput[tl](13.22,3.93){$t=\frac{k_n+1}{[n-(k_n+1)]}$}
\psline(-1,-8)(-0.74,-7.2)
\psline(7,1.03)(7.52,2.97)
\psline(7.67,1.03)(7.52,2.97)
\psline(9,1.03)(7.52,2.97)
\psline(11.5,1.03)(10.99,2.97)
\begin{scriptsize}
\psdots[dotstyle=*,linecolor=blue](5.5,-5.02)
\psdots[dotstyle=*,linecolor=blue](5.48,-4.46)
\psdots[dotstyle=*,linecolor=blue](5.48,-3.96)
\psdots[dotstyle=*,linecolor=blue](5.52,-1)
\psdots[dotstyle=*,linecolor=blue](5.5,-0.49)
\psdots[dotstyle=*,linecolor=blue](5.48,0)
\psdots[dotstyle=*,linecolor=darkgray](-1,-8)
\psdots[dotstyle=*,linecolor=darkgray](0.4,-8)
\psdots[dotstyle=*,linecolor=darkgray](3.48,-8)
\psdots[dotstyle=*,linecolor=darkgray](4.3,-8)
\psdots[dotstyle=*,linecolor=darkgray](5.44,-8)
\psdots[dotstyle=*,linecolor=darkgray](8.31,-8)
\psdots[dotstyle=*,linecolor=darkgray](12,-8)
\psdots[dotstyle=*,linecolor=darkgray](-0.74,-7.2)
\psdots[dotstyle=*,linecolor=darkgray](0.33,-7.2)
\psdots[dotstyle=*,linecolor=darkgray](1.2,-7.2)
\psdots[dotstyle=*,linecolor=darkgray](3.99,-7.2)
\psdots[dotstyle=*,linecolor=darkgray](7.5,-7.2)
\psdots[dotstyle=*,linecolor=darkgray](9.39,-7.2)
\psdots[dotstyle=*,linecolor=darkgray](0,-6)
\psdots[dotstyle=*,linecolor=darkgray](3.7,-6)
\psdots[dotstyle=*,linecolor=darkgray](8.15,-6)
\psdots[dotstyle=*,linecolor=darkgray](9,-6)
\psdots[dotstyle=*,linecolor=darkgray](-1.42,-3)
\psdots[dotstyle=*,linecolor=darkgray](-0.47,-3)
\psdots[dotstyle=*,linecolor=darkgray](0.42,-3)
\psdots[dotstyle=*,linecolor=darkgray](4.73,-3)
\psdots[dotstyle=*,linecolor=darkgray](5.48,-3)
\psdots[dotstyle=*,linecolor=darkgray](7.65,-3)
\psdots[dotstyle=*,linecolor=darkgray](8.25,-3)
\psdots[dotstyle=*,linecolor=darkgray](12.24,-3)
\psdots[dotstyle=*,linecolor=darkgray](-0.14,-1.38)
\psdots[dotstyle=*,linecolor=darkgray](6.02,-1.38)
\psdots[dotstyle=*,linecolor=darkgray](6.99,-1.38)
\psdots[dotstyle=*,linecolor=darkgray](10.29,-1.38)
\psdots[dotstyle=*,linecolor=darkgray](-1.27,1.03)
\psdots[dotstyle=*,linecolor=darkgray](0,1.03)
\psdots[dotstyle=*,linecolor=darkgray](0.64,1.03)
\psdots[dotstyle=*,linecolor=darkgray](3.16,1.03)
\psdots[dotstyle=*,linecolor=darkgray](7,1.03)
\psdots[dotstyle=*,linecolor=darkgray](7.67,1.03)
\psdots[dotstyle=*,linecolor=darkgray](9,1.03)
\psdots[dotstyle=*,linecolor=darkgray](11.5,1.03)
\psdots[dotstyle=*,linecolor=darkgray](-0.5,2.97)
\psdots[dotstyle=*,linecolor=darkgray](2.41,2.97)
\psdots[dotstyle=*,linecolor=darkgray](7.52,2.97)
\psdots[dotstyle=*,linecolor=darkgray](10.99,2.97)
\end{scriptsize}
\end{pspicture*}
\caption{Alternative model}
\label{fig:representationAM}
\end{figure}
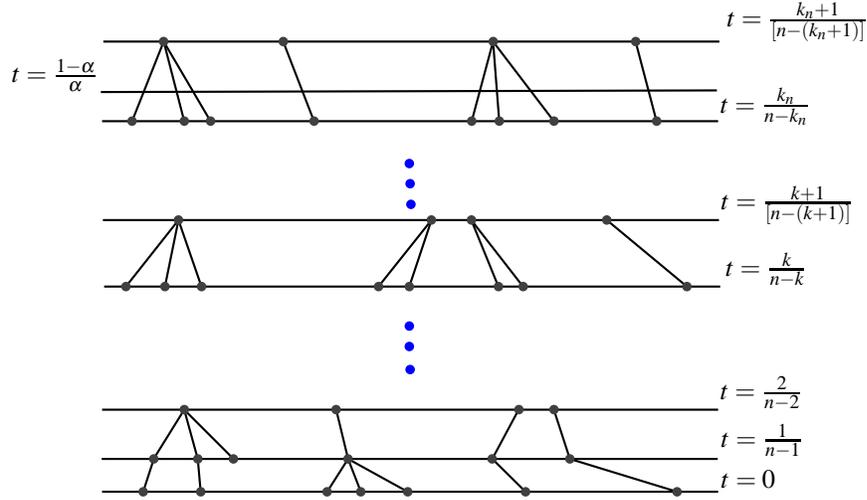

\subsection{Estimates for the tail of coalescence times}\label{calor}

Assume that we have two coalescing random paths starting at the same time at a (deterministic) distance apart. Denote by $\tau$ the meeting time of these two random paths. 
In order to verify conditions $(I), (B_1)$ and $(B_2)$ it is essential to have an estimate for the tail of its distribution of the following kind,

\begin{equation}
\mathbb{P}[\tau>t]\leq \frac{C}{\sqrt{t}}, \hspace{.3cm} t>0.\label{coalescingtime}
\end{equation}
The macroscopic random paths make jumps at times in 
\[
L_n=\left\{ 0, \frac{1}{n-1}, \frac{2}{n-2},...,\frac{k_n}{n-k_n},\frac{k_n+1}{n-k_n-1} \right\}.
\]
Under the diffusive scaling the microscopic trajectories may coalesce at times in 
$$
\hat{L}_n=\left\{ 0, \frac{n}{n-1}, \frac{2n}{n-2},...,\frac{nk_n}{n-k_n}, \frac{n(k_n+1)}{n-k_n-1} \right\}.
$$
Since

\begin{equation*}
1 \leq  \frac{n(j+1)}{n-j-1}- \frac{nj}{n-j}\leq C(\alpha) \ \mbox{for} \ j=0,1,2,...,k_n, 
\end{equation*}
where $C\left(\alpha\right)$ is a constant depending only on $\alpha$, we have that the distance between any two successive jump times is bounded from below and above. Therefore it suffices to study the case in which the jumps occur at non-negative integer numbers.

Now we formally define the coalescing time and prove inequality (\ref{coalescingtime}). To do so it is enough
to show the following case. Let $\hat{Q}=(\hat{Q}_k)_{k=0}^{\infty}$ be a collection of homogeneous, independent 
Poisson processes on $\mathbb{R}$. The processes $\hat{Q}_k, k=0,1, \ldots , k_n$, have intensities $(n-k)/n$ respectively, 
and the processes $\hat{Q}_k$ with $k > k_n$ have intensities $(n-k_n)/n$. Consider the subset of $\mathbb{R}^2$ given by 
$\cup_{k=0}^{\infty}\hat{Q}_k\times\{ k \}$.
For $u=(u_0,k)\in \mathbb{R}\times\{ k \}$, let $\tilde{u}^{1,k}={\arg\min}_{v\in \hat{Q}_{k+1}\times\{k+1 \}}\|u-v \|$. 
Define, inductively,

\[
\tilde{u}^{j,k}={\arg\min}_{v\in \hat{Q}_{k+j+1}\times\{ k+j+1 \}}\| \tilde{u}^{j-1,k}-v \|\quad j=1,2,...
\]
where $\tilde{u}^{0,k}=u$ . Let $X^{u}$ be the macroscopic trajectory obtained by linearly interpolating 
the points $\{u,\tilde{u}^{1,k},...,\tilde{u}^{j,k},... \}$.

Let $u=(u_1,0)$, $w=(w_1,0) \in\mathbb{R}\times\{ 0\}$ be such that $u_1<w_1$. Consider the trajectories 
$X^{u}$ and $X^{w}$, and, for $t\geq 0$, define

\begin{equation}
Z_t=Z_t(u,w)=X_t^{u}-X_t^{w}. \nonumber
\end{equation}
Assume that $w_1-u_1=1$ and define

\begin{equation}
\tau:=\min\{t\geq 0 : Z_t=0   \}. \nonumber
\end{equation}

\begin{prop}\label{tempode}
There exists a positive constant $C$ such that $\mathbb{P}[\tau >t]\leq \frac{C}{\sqrt{t}}$.
\end{prop}

\begin{proof} 
Since $\{Z_t(u,v),t\geq 0\}$ is a non-negative Martingale in $L^2$, Skorohod representation holds (see \cite{D}, page $384$ ). 
Then, there exists a Brownian motion with coefficient of diffusion $1$ starting from $1$ and stopping times $0=T_0,T_1,T_2,\ldots $ satisfying

\begin{equation}
Z_t\overset{\mathcal{D}}{=}B[T_t] \nonumber
\end{equation}
where $0=T_0,T_1,T_2,...$ are such that

\begin{equation}\label{a2}
T_t=\inf\{ s\geq T_{t-1} : B[s]-B[T_{t-1}]\notin \left( U_t(B(T_{t-1})), V_t(B[T_{t-1}])   \right)  \}, 
\end{equation}
where $\{ (U_t(r),V_t(r)),t\geq 1, r\in [0,\infty) \}$ is a family of random independent vectors. Note that for any $r\geq 0$, $t,t^{\prime}\geq k_n$, $(U_t(r),V_t(r))\overset{\mathcal{D}}{=}(U_{t^{\prime}}(r),V_{t^{\prime}}(r))$ and
$(U_t(r),V_t(r))\in [-r,0)\times (0,\infty)$ a.s.

Let $\tau^{\prime}:=\inf\{t\geq 0 : B[t]=0  \}$. Note that $\tau>t$ if and only if $\tau^{\prime}>T_t$. Then,

\begin{equation}\label{a1}
\mathbb{P}[\tau>t]=\mathbb{P}[\tau^{\prime}>T_t].
\end{equation}
Let $\varsigma > 0$ be a constant which will be specified later. Then, exactly as in \cite{CFD}, we get

\begin{eqnarray}\label{f1}
\mathbb{P}[\tau^{\prime}>T_t]& =& \mathbb{P}[\{\tau^{\prime}> T_t\}\cap \{T_t>\varsigma t \}]+\mathbb{P}[\{\tau^{\prime}>T_t\}\cap \{T_t\leq \varsigma t \}]\nonumber\\
                        & \leq & \mathbb{P}[\tau^{\prime}>\varsigma t]+\mathbb{P}[\{\tau^{\prime}>T_t\}\cap \{T_t\leq \varsigma t \}]\nonumber\\
                        & \leq & \frac{c_0}{\sqrt{t}}+\mathbb{P}[\{\tau^{\prime}>T_t\}\cap \{T_t\leq \varsigma t \}]\nonumber\\
\end{eqnarray}
where $c_0=c_0(\varsigma)\in(0,\infty)$. For $\lambda > 0$, we have

\begin{eqnarray}\label{relax}
\mathbb{P}[\{\tau^{\prime}>T_t\}\cap \{T_t\leq \varsigma t \}] & \leq & \mathbb{P}[\{\tau^{\prime}>T_t\}\cap \{e^{-\lambda T_t}\geq e^{-\lambda\varsigma t }\}]\nonumber\\
                                                         & = & \mathbb{P}[\{\Pi_{i=1}^t1_{\{Z_i>0 \}}>0\}\cap \{e^{-\lambda T_t}\geq e^{-\lambda\varsigma t}\}]\nonumber\\
                                                         & = & \mathbb{P}[\Pi_{i=1}^t1_{\{Z_i>0 \}}e^{-\lambda T_t}\geq e^{-\lambda\varsigma t}]\nonumber\\
                                                         & \leq & e^{\lambda\varsigma t}\mathbb{E}[\Pi_{i=1}^t 1_{\{Z_i>0 \}}e^{-\lambda T_t }].
\end{eqnarray}
Note that $T_t=\sum_{i=1}^t(T_i-T_{i-1})=\sum_{i=1}^tS_i(Z_{i-1})$ where $\{S_i(r),i\in\mathbb{N},r>0\}$ are independent random variables. However, for fixed $r > 0$, the random variables $\{S_i(r), i\in\mathbb{N}\}$ are not identically distributed. Then, $\mathbb{E}[\Pi_{i=1}^t1_{\{Z_i>0 \}}e^{-\lambda T_t}]$ equals

\begin{equation}
\mathbb{E}\left[ \mathbb{E}\left( \left.\Pi_{i=1}^{t-1}1_{\{Z_i>0 \}}\exp\left( -\lambda \sum_{i=1}^{t-1}S_i(Z_{i-1}) \right)1_{\{Z_t>0  \}}\exp\left( -\lambda S_t(Z_{t-1}) \right)\right| \mathcal{F}_{t-1} \right)\right]\nonumber
\end{equation}
which is bounded from above by

\begin{equation}
\label{relax1}
\Pi_{i=1}^t \sup_{r>0}\mathbb{E}\left[ 1_{\{Z_i>0\}}e^{-\lambda S_i(r) }\right],
\end{equation}
where $\left\{ \mathcal{F}_t \right\}$ is the $\sigma$-algebra generated by the random variables $Z_0,Z_1,...,Z_t$.

Since $\mathbb{E}\left[ 1_{\{Z_t>0\}}e^{-\lambda S_{t}(r)}\right]\leq \mathbb{P}_r[Z_t>0]$ where $\mathbb{P}_r[Z_t>0]$ is the probability that two paths starting from distance $r$ in the level $t-1$, have not coalesce by level $t$. Then

\begin{eqnarray}\label{4j}
\sup_{r\leq 10}\mathbb{E}\left[ 1_{\{Z_t>0\}}e^{-\lambda S_1(r) } \right]&\leq& \sup_{r\leq 10} \mathbb{P}[Z_t>0]\nonumber\\
                  &=&\mathbb{P}_{10}[Z_t>0].\nonumber\\
\end{eqnarray}
The previous equality follows from the translation invariance of the Poisson processes and the fact that $r \leq 10$.

For any fixed $x\in\mathbb{R}$ define $p(t,x):=\arg\min_{y\in \hat{Q}_t}|y-x|$. Let $B_t:=\{ p(t,0)\neq p(t,10) \}$. By the translation invariance of the Poisson processes we have that

\[
\mathbb{P}[B_t]=\mathbb{P}_{10}[Z_t>0 ].
\]
Consider the event

\[
C_t=\{|[-10,0]\cap \hat{Q}_t|=0\}\cap\{|(0,10)\cap \hat{Q}_t|=1\}\cap \{|[10,20]\cap \hat{Q}_t|=0\}.
\]
Then, $C_t\subset B_t^c$. Let $\kappa_t$ be the intensity of the Poisson process $\hat{Q}_t$. Since $\kappa_t \in [\alpha,1]$, we have that

\begin{equation*}
\mathbb{P}[C_t]=e^{-20\kappa_t }e^{-10\kappa_t}\kappa_t 10\geq e^{-30}\alpha 10=c_0>0.
\end{equation*}
Then, $\mathbb{P}_{10}[Z_t>0 ]\leq 1-c_0=c_1<1$. It follows from (\ref{4j}) that,

\begin{equation}\label{descansar}
\sup_{r\leq 10}\mathbb{E}\left[ 1_{\{Z_t>0\}}e^{-\lambda S_t(r) }\right]\leq c_1<1.
\end{equation}
where $c_1$ does not depend on $t$. Now, we will show that there exists a positive constant $c_{2}$ independent of $t$ such that

\begin{equation}\label{pao}
\sup_{r\geq 10}\mathbb{E}\left[e^{-\lambda S_t(r) }\right]\leq c_2.
\end{equation}
Unless otherwise stated, from now $r\geq 10$. The strategy to show (\ref{pao}) is to find a suitable set $R_{\epsilon_0}$ such that

\begin{equation}
\mathbb{E}\left[ e^{-\lambda S_t(r)} \right]\leq \mathbb{P}[(U_t(r),V_t(r))\in R_{\epsilon_0}](1-a(\epsilon_0))+a(\epsilon_0), \nonumber
\end{equation}
where $a(\epsilon_0)$ and $\mathbb{P}[(U_t(r),V_t(r))\in R_{\epsilon_0}]\in (0,1) $ do not depend on $t$ and $r$.

Let $\epsilon > 0$ be given. Define the sets $A_\epsilon:=(-\epsilon,0)\times(0,\infty)$ and $B_\epsilon:=[-r,0)\times(0,\epsilon)$. Also, 
let $F^t(x)=\mathbb{P}_r[Z_t-r\leq x]$. Then (see \cite{D} page 403),

\begin{eqnarray}
\label{4a}
\mathbb{P}[(U_t(r),V_t(r))\in A_\epsilon\cup B_\epsilon]  &=&  \frac{1}{c}\int\int_{(u,v)\in A_\epsilon\cup B_\epsilon}(v-u)dF^t(u)dF^t(v) \\ 
&\leq&  \frac{1}{c}\int\int_{(u,v)\in (-\epsilon,0)\times(0,\infty)}(v-u)dF^t(u)dF^t(v) \nonumber \\
&+& \frac{1}{c}\int\int_{(u,v)\in[-r,0)\times(0,\epsilon)}(v-u)dF^t(u)dF^t(v) \nonumber
\end{eqnarray}
where $c=c(r)=\int_{-r}^{0}(-u)dF^t(u)=\int_{0}^{\infty}vdF^t(v)$. By the properties of the Poisson process we have $c<\infty$. A straightforward computation yields that

\begin{equation}
\frac{1}{c} \int\int_{(u,v)\in (-\epsilon,0)\times(0,\infty)} (v-u) dF^t(u)dF^t(v) \nonumber
\end{equation}
is less than or equal to
%\\
%= \frac{1}{c}\int\int_{(u,v)\in (-\epsilon,0)\times(0,\infty)}vdF^t(u)dF^t(v) 
%\\
%+ \frac{1}{c}\int\int_{(u,v)\in (-\epsilon,0)\times(0,\infty)}(-u)dF^t(u)dF^t(v)  
%\\ 
\begin{equation}
\leq \int_{u\in(-\epsilon,0)}dF^t(u)+\frac{1}{c}\int_{u\in(-\epsilon,0)}(-u)dF^t(u)
\end{equation}
which equals
\begin{equation}
\mathbb{P}_r[-\epsilon\leq Z_t-r\leq 0]+\frac{1}{c}\int_{u\in(-\epsilon,0)}(-u)dF^t(u). \nonumber
\end{equation}
Therefore,

\begin{equation}\label{4b}
\frac{1}{c}\int\int_{(u,v)\in (-\epsilon,0)\times(0,\infty)}(v-u)dF^t(u)dF^t(v)\leq \mathbb{P}_r[-\epsilon\leq Z_t-r\leq 0]+\frac{\epsilon}{c}.
\end{equation}
The same reasoning yields

\begin{equation}\label{4c}
\frac{1}{c}\int\int_{(u,v)\in[-r,0)\times(0,\epsilon)}(v-u)dF^t(u)dF^t(v)\leq \mathbb{P}_r[0\leq Z_t-r\leq \epsilon]+\frac{\epsilon}{c}.
\end{equation}
Collecting (\ref{4a}), (\ref{4b}) and (\ref{4c}) we get

\begin{equation}\label{4d}
\mathbb{P}[(U_t(r),V_t(r))\in A_\epsilon\cup B_\epsilon]\leq \mathbb{P}_r[  -\epsilon\leq Z_t-r\leq \epsilon]+\frac{2\epsilon}{c}.
\end{equation}
Again, by the translation invariance of the Poisson processes, we get
\begin{equation}\label{3.69}
\mathbb{P}_r[-\epsilon\leq Z_t-r\leq \epsilon]=\mathbb{P}[-\epsilon\leq p(t,r)-p(t,0)-r\leq \epsilon].
\end{equation}
Now define the events $C_1(t):=\{ [-10/4,10/4]\cap \hat{Q}_t\neq \emptyset \}$ and $C_2(t):=\{ [-10/4+r,r+10/4]\cap \hat{Q}_t
\neq \emptyset \}$. Then,

\begin{multline}
\mathbb{P}[-\epsilon\leq p(t,r)-p(t,0)-r\leq \epsilon] \leq \\  
\mathbb{P}[-\epsilon\leq p(t,r)-p(t,0)-r\leq \epsilon| C_1\cap C_2]\mathbb{P}[C_1\cap C_2] + \mathbb{P}[{C_1}^c\cup {C_2}^c]\nonumber
\end{multline}
which equals
\begin{equation*}
\mathbb{P}[-\epsilon\leq p(t,r)-p(t,0)-r\leq \epsilon| C_1\cap C_2](1-e^{-5\kappa_t })^2 + [1-(1-e^{-5\kappa_t })^2]. 
\end{equation*}
Now note that, conditioned on $[C_1\cap C_2], p(t,r)$ and $p(t,0)$ are independent. Therefore,

\begin{equation*}
\mathbb{P}[-\epsilon\leq p(t,r)-p(t,0)-r\leq \epsilon| C_1\cap C_2]\leq \mathbb{P}[|[-\epsilon/2,\epsilon/2]\cap \hat{Q}_t|\geq 1].
\end{equation*}
Thus,
\begin{equation}\label{nana}
\mathbb{P}[-\epsilon\leq p(t,r)-p(t,0)-r\leq \epsilon| C_1\cap C_2]\leq 1-e^{-\kappa_t\epsilon}.
\end{equation}
It follows from (\ref{4d}),(\ref{3.69}) and (\ref{nana}) that
\begin{equation*}
\sup_{r\geq 10}\mathbb{P}[(U_t(r),V_t(r))\in A_\epsilon\cup B_\epsilon]\leq (1-e^{-\kappa_t\epsilon})(1-e^{-5\kappa_t })^2+
[1-(1-e^{-\kappa_t })^2]+\frac{2\epsilon}{c}.
\end{equation*}
Since $\kappa_t\in \left[ \alpha,1 \right]$, we get

\begin{equation}\label{c(r)}
\sup_{r\geq 10}\mathbb{P}[(U_t(r),V_t(r))\in A_\epsilon\cup B_\epsilon]\leq (1-e^{-\epsilon})(1-e^{-5})^2+[1-(1-e^{-5\alpha })^2]+\frac{2\epsilon}{c}.
\end{equation}
Note that $c=c(r)=\int_{-r}^{0}(-u)dF^t(u)=\int_{0}^{\infty}vdF^t(v)$ depends on both, $t$ and $r$. However, for $r \geq 10$ (in fact $10$ is just to fix ideas), $c$ can be bounded from below by constant not depending on $t$, neither on $r$. Indeed,

\begin{eqnarray}
c=\int_{0}^\infty vdF^t(v) & \geq &\int_1^\infty vdF^t(v)\nonumber\\
                           & \geq & \mathbb{P}_r[Z_t-r\geq 1].\nonumber\\
\end{eqnarray}
Again, by the translation invariance of the Poisson processes, we have

$$
\mathbb{P}_r[Z_t-r\geq 1]=\mathbb{P}[p(r)-p(0)\geq r+1].
$$
Let $D=\left\{ |[-2,-1]\cap \hat{Q}_t|=1 \right\}\cap \left\{ |[-1,10]\cap \hat{Q}_t|=0 \right\}\cap\left\{ |[r-1,r]\cap \hat{Q}_t|=1 \right\}$. Then, $D\subset \{p(r)-p(0)\geq r+1\}$. Therefore,

\begin{equation}\label{c(r)1}
\mathbb{P}_r[Z_t-r\geq 1]\geq \mathbb{P}[D]=e^{-11\kappa}(\kappa_t)^2e^{-2\kappa_t}=e^{-13\kappa_t}(\kappa_t)^2\geq e^{-13}\alpha^2.
\end{equation}
It follows from (\ref{c(r)})-(\ref{c(r)1}) that

\begin{equation*}
\sup_{r\geq 10}\mathbb{P}[(U_t(r),V_t(r))\in A_\epsilon\cup B_\epsilon]\leq (1-e^{-\epsilon})(1-e^{-5})^2+[1-(1-e^{-5\alpha })^2]+\frac{2\epsilon e^{13}}{\alpha^2}.
\end{equation*}
Let $\epsilon_0$ be such that

$$
(1-e^{-\epsilon_0})(1-e^{-5})^2+[1-(1-e^{-5\alpha })^2]+\frac{2\epsilon_0 e^{13}}{\alpha^2}=c_2(\epsilon_0)<1 .
$$
Thus,

\begin{equation}\label{casi}
\sup_{r\geq 10}\mathbb{P}[(U_t(r),V_t(r))\in A_{\epsilon_0}\cup B_{\epsilon_0}]\leq c_2(\epsilon_0).
\end{equation}
Now, for $r\geq 10$, we have that

\begin{eqnarray}
\mathbb{E}[1_{\{ Z_t>0 \}}e^{-\lambda S_t(r)}]  
%&\leq&  \mathbb{E}\left[ e^{-\lambda S_t(r)}\right] \nonumber \\
%&\leq& \mathbb{E}\left[ e^{-\lambda S_t(r)}1_{\{(U_t,V_t)\in A_{\epsilon_0}\cup B_{\epsilon_0} \}}\right] \nonumber \\
%&+& \mathbb{E}\left[ e^{-\lambda S_t(r)}1_{\{(U_t,V_t)\in (A_{\epsilon_0}\cup B_{\epsilon_0})^c \}}\right] \nonumber \\
%&=&
&\leq& \mathbb{E}\left[ \left.\mathbb{E}\left( e^{-\lambda S_t(r)}1_{\{(U_t,V_t)\in A_{\epsilon_0}\cup B_{\epsilon_0} \}} \right| (U_t(r),V_t(r)) \right) \right] \nonumber \\
&+&\mathbb{E}\left[ \left.\mathbb{E}\left( e^{-\lambda S_t(r)}1_{\{(U_t,V_t)\in (A_{\epsilon_0}\cup B_{\epsilon_0})^c \}} \right| (U_t(r),V_t(r)) \right)\right] \nonumber \\
&\leq& \mathbb{P}[(U_t(r),V_t(r))\in A_{\epsilon_0}\cup B_{\epsilon_0}]\nonumber \\
&+&\mathbb{E}\left[ \left.\mathbb{E}\left( e^{-\lambda S_t(r)}1_{\{(U_t,V_t)\in [-r,-\epsilon_0)\times [\epsilon_0,\infty)\}} \right| (U_t(r),V_t(r)) \right)\right]. \nonumber
\end{eqnarray}
Let $(U_t,V_t)$ be given. Also, let $S_t(r)$ be the stopping time of $B$ determined by $(U_t(r), \newline V_t(r))$. Then, if $T_{\epsilon_0}:=\inf\{t>0 : B[t]\notin (-\epsilon_0, \epsilon_0)  \}$ we have (see \cite{D} page 400) that

\small 
\begin{eqnarray}
\mathbb{E}[e^{-\lambda S_t(r)}] &\leq&  \mathbb{P}[(U_t(r),V_t(r))\in A_{\epsilon_0}\cup B_{\epsilon_0}]+ \mathbb{E}[e^{-\lambda T_{\epsilon_0}}]
\mathbb{P}[(U_t(r),V_t(r))\in (A_{\epsilon_0}\cup B_{\epsilon_0})^c] \nonumber \\
%&\leq& \mathbb{P}[(U_t(r),V_t(r))\in A_{\epsilon_0}\cup B_{\epsilon_0}]+\frac{1}{\cosh(\sqrt{2\lambda}\epsilon_0)}[1-\mathbb{P}[(U_t(r),V_t(r))\in 
%A_{\epsilon_0}\cup B_{\epsilon_0}]] \nonumber \\
%&=&
&\leq& \mathbb{P}[(U_t(r),V_t(r))\in A_{\epsilon_0}\cup B_{\epsilon_0}]\left(1-\frac{1}{\cosh(\sqrt{2\lambda}{\epsilon_0})}  \right)+
\frac{1}{\cosh(\sqrt{2\lambda}{\epsilon_0})}. \nonumber
\end{eqnarray}
\normalsize
It follows from (\ref{casi}) that
 
\begin{eqnarray}\label{4h1}
\sup_{r\geq 10}\mathbb{E}[e^{-\lambda S_t(r)}] & \leq & \sup_{r\geq 10}\mathbb{P}[(U_t(r),V_t(r))\in A_{\epsilon_0}\cup B_{\epsilon_0}]\left(1-\frac{1}{\cosh(\sqrt{2\lambda}{\epsilon_0})}\right) \nonumber \\
&+& \frac{1}{\cosh(\sqrt{2\lambda}{\epsilon_0})}\nonumber\\
& \leq & c_1(\epsilon_0)\left(1-\frac{1}{\cosh(\sqrt{2\lambda}{\epsilon_0})}  \right)+\frac{1}{\cosh(\sqrt{2\lambda}{\epsilon_0})} \nonumber\\
\end{eqnarray}
Let $c_3:=c_3(\epsilon_0)=c_2(\epsilon_0)\left(1-\frac{1}{\cosh(\sqrt{2\lambda}{\epsilon_0})}  \right)+\frac{1}{\cosh(\sqrt{2\lambda}{\epsilon_0})}<1$ and let $c_4=c_4(\epsilon_0) \newline =c_3(\epsilon_0)\vee c_1<1$. Then, from (\ref{4h1}) and (\ref{descansar}) we get

\begin{equation}
\sup_{r\leq 10}\mathbb{E}\left[1_{\{ Z_t>0 \}} e^{-\lambda S_t(Z_{t-1})}|Z_{t-1}=r \right]\leq c_1\quad \nonumber
\end{equation}
and

\begin{equation}
\sup_{r\geq 10}\mathbb{E}\left[1_{\{ Z_t>0 \}} e^{-\lambda S_t(Z_{t-1})}|Z_{t-1}=r \right]\leq c_3. \nonumber
\end{equation}
Therefore,

\begin{equation*}
\sup_{r>0}\mathbb{E}\left[1_{\{Z_t>0\}} e^{-\lambda S_t(Z_{t-1})}|Z_{t-1}=r \right]\leq c_4<1.
\end{equation*}
Note that $c_4\in(0,1)$ depends only on $\epsilon_0$. Thus,

\begin{equation}
\max_{1\leq i\leq t}\sup_{r>0}\mathbb{E}\left[1_{\{Z_i>0\}} e^{-\lambda S_i(Z_{i-1})}|Z_{i-1}=r \right]<c_4. \nonumber
\end{equation}
It follows from (\ref{relax}) and (\ref{relax1}) that

\begin{equation*}
\mathbb{P}[\tau >t] \leq \frac{c_0(\varsigma)}{\sqrt{t}}+\left[ e^{\lambda\varsigma}c_4 \right]^t.
\end{equation*}
Let $\varsigma_0>0$ be such that $ e^{\lambda\varsigma_0}c_4=c_5<1$. Note that $c_5$ depends only on $\epsilon_0$ and $\varsigma_0$. Therefore,

\begin{equation*}
\mathbb{P}[\tau >t] \leq \frac{c_0(\varsigma_0)}{\sqrt{t}}+(c_5(\epsilon_0,\varsigma_0))^t.
\end{equation*}
Since $c_5\in(0,1)$, there exists $c_6=c_6(\epsilon_0,\varsigma_0)$ such that $(c_5)^t\leq \frac{c_6}{\sqrt{t}}$. Then, taking $c_7=c_7(\epsilon_0,\varsigma_0)=c_0 \vee c_6$ we get

\begin{equation}
\mathbb{P}[\tau > t]\leq \frac{c_7}{\sqrt{t}}. \nonumber
\end{equation}
\end{proof}

\begin{lem}\label{tempo} Let $k \in \mathbb{N}$. Then, for every $u,v\in\mathbb{R}\times\{ k\}$, 
such that $v_1-u_1=1$ let $Z_t:=Z_t(u,v)=X_{k+t}^v-X_{t+k}^{u}$ and
$\tau_{k}:=\min\{t>k : Z_t=0 \}$. Then, there exists a positive constant $C$ such that

\[
\mathbb{P}[\tau_k>t]\leq \frac{C}{\sqrt{t}}.
\]
\end{lem}

\begin{proof} From the previous proposition we know that $\mathbb{P}[\tau_0>t]\leq \frac{C}{\sqrt{t}}$
where $C$ is a positive constant independent of levels $0$ and $t$. The Lemma follows from the previous 
theorem and the fact that the intensities belong to $[\alpha, 1]$.
\end{proof}

\begin{lem} \label{tempo2} Let $u=(x,0)$ and $v=(y,0)$ with $x<y$. Then,

\begin{equation}
\mathbb{P}[Z_t(u,v)>0]\leq \frac{C(y-x)}{\sqrt{t}}. \nonumber
\end{equation}
\end{lem}

\begin{proof} The following proof is analogous to that of Lemma $3.2$, page 31 in \cite{SUN}.
 If $y-x\leq 1$, the result follows from the previous lemma. Now, let $y-x > 1 $. Since

\begin{equation}
\left\{ Z_t(u,v)>0 \right\}\subset \left\{ \cup_{i=0}^{\lceil y\rceil-\lfloor x\rfloor-1}\left\{ Z_t((\lfloor x\rfloor + i,0),(\lfloor x\rfloor+i+1)) >0 \right\}
\right\} \nonumber
\end{equation}
we get

\begin{eqnarray*}
\mathbb{P}[ Z_t(u,v)>0] & \leq & \sum_{i=0}^{\lceil y\rceil-\lfloor x\rfloor-1}\mathbb{P}[Z_t((\lfloor x\rfloor + i,0),(\lfloor x\rfloor+i+1)) >0]\\
            &  =&(\lceil y\rceil-\lfloor x\rfloor-1)\mathbb{P}[Z_t((0,0),(1,0)) >0]\\
             & \leq &2(y-x)\mathbb{P}[Z_t((0,0),(1,0)) >0],\\
\end{eqnarray*}
where the last equality follows form the translation invariance of Poisson processes. 
Then, the results follows from the previous lemma. 
\end{proof}

\subsection{Verification of condition $B_2$} \label{vcb2}

 In order to prove Condition $B_2$ we need to estimate the number of points of the point process in an interval of length 
$\epsilon$. This is the content of the following Lemma.

\begin{lem}
\label{cansado} 
Let $k=1,2,...,k_n$ and make $A_k:=\left\{ |\hat{Q}_k\cap [0,\epsilon\sqrt{n}]|\in\left[\frac{\alpha}{2}\epsilon\sqrt{n},
\frac{2e}{\alpha}\epsilon\sqrt{n} \right]   \right\}$ for some $\epsilon>0$ fixed, that is $A_k$ is
the event that the number of points of the Poisson processes on $[0,\epsilon\sqrt{n}]$ belongs to $\left[\frac{\alpha}{2}\epsilon\sqrt{n},
\frac{2e}{\alpha}\epsilon\sqrt{n} \right]$. Then 
\begin{equation}
\limsup_{n\rightarrow\infty}\sup_{1\leq k\leq k_n}\mathbb{P}[A_k^c]=0. \nonumber
\end{equation}
Thus, the probability that $\cup_{k=1}^{k_n}A_k^c$ occurs infinitely often is 0.

\end{lem}

\begin{proof} Since $(T_i)_{i=1}^{\left\lceil \frac{\alpha}{2}\epsilon\sqrt{n} \right\rceil}$ is a collection 
of exponential i.i.d. random variable with rate $\frac{n-k}{n}\in[\alpha,1]$ we have

\[
\mathbb{P}\left[ |\hat{Q}_k\cap [0,\epsilon\sqrt{n}]|<\frac{\alpha}{2}\epsilon\sqrt{n} \right]
\leq \mathbb{P}\left[  \sum_{i=1}^{\left\lceil \frac{\alpha}{2}\epsilon\sqrt{n} \right\rceil}T_i>\epsilon\sqrt{n} \right].
\]
A standard large deviation argument gives exponential decay with $n$ of this probability. Also, we have

\[
\mathbb{P}[|\hat{Q}_k\cap [0,\epsilon\sqrt{n}]|>\frac{2e}{\alpha}\epsilon\sqrt{n}]\leq 
\mathbb{P}\left[  \sum_{i=1}^{\left\lceil \frac{2e}{\alpha}\epsilon\sqrt{n} \right\rceil}T_i<\epsilon\sqrt{n} \right].
\]
Since $M_n=\sum_{i=1}^{\left\lceil \frac{2e}{\alpha}\epsilon\sqrt{n} \right\rceil}T_i$ is 
Gamma distributed with parameters $\left\lceil \frac{2e}{\alpha}\epsilon\sqrt{n} \right\rceil$ and $\frac{n}{n-k}$, 
it follows from the Stirling approximation 
that $\mathbb{P}[M_n<\epsilon\sqrt{n}]$ converges exponentially fast to zero.

By the above, we conclude that $\mathbb{P}[A_k^c]$ decays exponentially fast with $n$. 
Since this is true for every $k=0,1,...,k_n$, the proof is complete. 
\end{proof} 

The lemma above implies that with high probability the number of macroscopic trajectories crossing a interval of length $\epsilon\sqrt{n}$ are of order $\sqrt{n}$.    

Now, we establish an FKG inequality which is the main tool used in the proof of condition B2. The proof that this inequality holds proceeds in two steps. In the first step we introduce a random model whose paths make jumps in a discrete spatial set and where the FKG inequality holds. In the second step, we show that the desired FKG inequality holds by taking an appropiate limit. Let $r > 0$ be fixed and consider the set $\mathbb{Z}_r:=\{ rk : k\in\mathbb{Z}\}$. For $k=0,1,2,...,k_n$, let $E_k := E_k(r)=\mathbb{Z}_r\times \{ l_k \}$, $l_k=\frac{nk}{n-k}$ and $E := \cup_{k=0}^{k_n}E_k$. Say that $z_k=(rj,l_k)\in E_k, j\in\mathbb{Z}$ is open if $Q_k\cap [rj,r(j+1))\neq \emptyset$ and is closed in any other case. If $p_k^r$ denotes the probability of $z_k$ being open, then 

\begin{equation}
p_k^r = \mathbb{P}[Q_k\cap [rj,r(j+1))\neq \emptyset] = 1-e^{-\frac{n-k}{n}r}. \nonumber
\end{equation}

Let $\Omega^{n,r}=\left(\omega(z), z\in E  \right)$ be a family of independent, Bernoulli random variables. For $k=0,1,2,...,k_n$, assume that the parameter of the Bernoulli variables at lever $l_k$ is $p_k^r$. Also, consider a second family of independent, Bernoulli random variables $\Gamma=(v(z_k),z_k\in E_k)_{0\leq k\leq k_n-1}$, each one with parameter $1/2$.

In what follows we consider an auxiliary system of coalescing paths itroduced by Gangopadhyay, Roy and Sarkar \cite{GRS}. This system is known as drainage network and its scaling limit was obtained by Coletti, Fontes and Dias \cite{CFD}. It evolves in a space-time lattice in the following way: Each space-time point $\left(x,t\right)$ of the lattice is consider open or closed according to an iid family of Bernoulli random variables; if at time $t$ a walk is at space site $x$ then at time $t+1$ its position is the open space site nearest to $x$ in the next level. By next level we mean a space-time site with time coordinate equal to $t+1$. We denote by $h\left((x,t)\right)$ this closest point. If we have two possible choises for the nearest open site, we choose one of them with probability $1/2$.

More  precisely, for $z_k=(a_k,l_k)\in E_k$, $(k<k_n)$, let $h_k(z_k)={\arg\min_{z^{\prime}\in E_{k+1}}}\|z-z^{\prime} \|$ and $w(h_k(z)) = 1$. In case of tie, proceed as follows. If $v(z_k)=1$, let $h_k(z_k)$ be defined as the nearest point to its right in the next level. If $v(z_k)=0$, let $h_k(z_k)$ be defined as the nearest point to its left in the next level. 
Let $h_k^{0}(z_k) := z_k$. For $k\leq j\leq k_n-k$, let

$$
h_k^j(z_k) := {\arg\min}_{z^{\prime}\in E_{k+j}}\|h_k^{(j-1)}(z_k)-z^{\prime} \| \quad \wedge \quad w(h_k^j(z_k))=1.
$$
Consider the oriented graph $\mathcal{G}_r=(V,\mathcal{E})$ with vertices set $V=E$ and edge set
$$
\mathcal{E}=\left\{ (u,h_k^{1}(u)): u\in E_k \right\}_{k=0}^{k_n}.
$$

For $a\in\mathbb{R}$, let $m_r(a)=r\left\lfloor \frac{a}{r} \right\rfloor$. Denote by $J_{a}^r$ the path in $\mathcal{G}_r$ starting from $(m_r(a),  0)$. Let 
$a\in[-n^{\delta/2+1/2},n^{\delta/2+1/2}]$ and 
$$
S^n_{a}:=\left\{(\sqrt{n}x_1,nx_2) : (x_1,x_2)\in \beta^n_{(a/\sqrt{n},0)}  \right\}.
$$
By abuse of notation we write $\beta^n_{(a/\sqrt{n},0)}$ even if, almost surely, $(a,0)\notin \hat{Q}_0$. Observe that the increments of $S_{a}^n$  depend on the Poisson 
point processes $\hat{Q}_k$ with rate $\frac{n-k}{n}$, $k=0,1,2,...,k_n$.

Let $t\in(0,1/\alpha-1)$ be fixed. Denote the state space of $J_{a}^r$ by $\tilde{\Pi}_{a}^r$. 
Introduce a partial order $\prec$ on $\tilde{\Pi}_{a}^r$ as follows. Let $\pi_1^r,\pi_2^r\in \tilde{\Pi}_{a}^r$ be given. Say that

\begin{equation}
\pi_1^r\prec\pi_2^r \Leftrightarrow \pi_1^r(l_k)-\pi_1^r(l_s)\leq \pi_2^r(l_k)-\pi_2^r(l_s), \nonumber
\end{equation}
for any $\left\lfloor nt\right\rfloor\geq k\geq s\geq 0$.

\begin{prop}\label{FKG}
Let $\pi_1^r,\pi_2^r\in\tilde{\Pi}_a^r$ be such that $\pi_1^r\prec\pi_2^r$. Let $x<a<y$. Then,

\begin{equation}\label{83}
\mathbb{P}[J_x^r(l_{\left\lfloor nt\right\rfloor})<J_a^r(l_{\left\lfloor nt\right\rfloor})| J_a^r=\pi_1^r]\leq \mathbb{P}[J_x^r(l_{\left\lfloor nt\right\rfloor})<J_a^r(l_{\left\lfloor nt\right\rfloor})| J_a^r=\pi_2^r]
\end{equation}
and

\begin{equation}\label{84}
\mathbb{P}[J_y^r(l_{\left\lfloor nt\right\rfloor})>J_a^r(l_{\left\lfloor nt\right\rfloor})| J_a^r=\pi_1^r]\geq \mathbb{P}[J_y^r(l_{\left\lfloor nt\right\rfloor})>J_a^r(l_{\left\lfloor nt\right\rfloor})| J_a^r=\pi_2^r].
\end{equation}
\end{prop}

\begin{proof}
If $m_r(x)=m_r(a)$, then $J_x^r=J_a^r$. Therefore, for almost every realization $\pi^r$ of $\tilde{\Pi}_a^r$, we have

\[
\mathbb{P}[J^r_x(l_{\left\lfloor nt\right\rfloor})<J_a^r(l_{\left\lfloor nt\right\rfloor})|J_a^r=\pi^r]=0.
\]

This gives (\ref{83}). If $m_r(x)<m_r(a)$, then the proof follows in the same lines of the analogous result proved in \cite{CFD}, page $1189$. The proof of (\ref{84}) is entirely analogous to the proof of (\ref{83}).
\end{proof}

\begin{rem}
Proposition \ref{FKG} is the analogous of Proposition (1) in \cite{CFD}. 
\end{rem}

\begin{lem}\label{minana}
Let $a\in[-n^{\delta/2+1/2},n^{\delta/2+1/2}]$ be fixed.
Then, for any $k=0,1,...,k_n$, we have

\begin{equation}
J_a^r(l_k)\overset{a.s}{\underset{r\rightarrow 0^+}{\longrightarrow}}S_a^n(l_k). \nonumber
\end{equation}
\end{lem}

\begin{proof}
Fix $n \in \mathbb{N}$. Since a Poisson process has (a.s.) no limit points, we may conclude that (a.s.) if $\omega$ is a realization of the $k_n$´s Poisson processes, then there exists a constant $r\left(\omega\right)$ depending only on $\omega$  such that if $M_k(\omega) := |S^n_a(l_k)(\omega)-J_a^r(l_k)(\omega)|$, then $\max_{1\leq k\leq k_n }M_k(\omega)=0$. Therefore \begin{equation*}
J_a^r(l_k)\overset{a.s}{\underset{r\rightarrow0^+}{\longrightarrow}}S^n_a(l_k),\quad k=0,1,...,k_n,
\end{equation*}
which is the desired conclusion.
\end{proof}

In order to verify condition $B_2$ of Theorem \ref{omega4} we must prove that, for any $t \in (0,1/\alpha-1)$,

\begin{equation}\label{fe2}
\epsilon^{-1}\limsup_{n\rightarrow\infty}\sup_{(a,t_0)\in \mathbb{R}\times[0,1/\alpha-1]}
\mathbb{P}[\eta_{\beta^n}(t_0,t;a,a+\epsilon)\geq 3]\underset{\epsilon\rightarrow 0^+}{\longrightarrow}0.
\end{equation}
Let $r_n(t):=\max\left\{ k: \frac{k}{n-k}< 1/\alpha -1-t  \right\}$ and $m_k=\frac{k}{n-k}$. By
translation invariance, (\ref{fe2}) is equivalent to show that

\begin{equation}\label{fe3}
\epsilon^{-1}\limsup_{n\rightarrow\infty}\sup_{0\leq k \leq r_n(t) }
\mathbb{P}[\eta_{\beta^n}(m_k,t;0,\epsilon)\geq 3]\underset{\epsilon\rightarrow 0^+}{\longrightarrow}0.
\end{equation}

If $S^n:=\left\{ (\sqrt{n}x, ny): (x,y)\in \beta^n \right\}$ we have

\[
\mathbb{P}[\eta_{\beta^n}(m_k,t;0,\epsilon)\geq 3]=\mathbb{P}[\eta_{S^n}(nm_k,nt;0,\sqrt{n}\epsilon)\geq 3].
\]
Let $A_0=\{|\hat{Q}_0\cap [0,\epsilon\sqrt{n}] |=m\in\left[\alpha/2\epsilon\sqrt{n},
\frac{2e}{\alpha}\epsilon\sqrt{n}\right]\}$ for each $k\geq 0$. Denote by $x_1,x_2,...,x_m$ the marks of the Poisson  process and let

\begin{equation*}
\eta_0^{\prime}:=\eta_0^{\prime}(x_1,x_2,...,x_m)=|\{ S_j^n(tn): 1\leq j \leq m \}|,
\end{equation*}
where $S_j^n$ denote the microscopic path that start in $x_j$. Then,

\begin{equation}\label{casi1}
\mathbb{P}[\eta_{S^n}(0,nt;0,\sqrt{n}\epsilon)\geq 3|A_0]=\mathbb{P}[\eta_0^{\prime}\geq 3].
\end{equation}
At this point we need to enlarge the set of starting points of the random paths under consideration. We do so by adding deterministic points to the set of points $x_i$'s. More precisely, for $i=0,1,2,...,m$, let $y_i := i\frac{\epsilon\sqrt{n}}{m}$. Since $\frac{\epsilon\sqrt{n}}{m}\in \left[\frac{\alpha}{2e}, \frac{2}{\alpha} \right]$, the distance
between these points can not be larger than $\frac{2}{\alpha}$. Let 

$$
\{z_0=0,z_1,z_2,...,z_{2m},z_{2m+1}=\epsilon\sqrt{n}\}
$$
be the union of the $x^{\prime}s$ and $y^{\prime}s$ points ordered from left to right. Let $\hat{S}_j^n$ be the path starting from the point $z_j$. If $\eta_0^{\prime}\geq 3$, there exists $j$ with $1\leq j\leq 2m$ such that 
$\hat{S}^n_{j-1}(tn)<\hat{S}^n_j(tn)<\hat{S}^n_{2m+1}(tn)$. Therefore, by Lemma (\ref{minana}), we have

\begin{eqnarray*}
\mathbb{P}[\eta_0^{\prime}\geq 3] &\leq &\sum_{j=1}^{2m}\mathbb{P}[\hat{S}^n_{j-1}(tn)<\hat{S}^n_{j}(tn)<\hat{S}^n_{2m+1}(tn)]\\
                          & =   &\lim_{r\rightarrow 0^+}\sum_{j=1}^{2m}\mathbb{P}[J_{j-1}^r(tn)<J_{j}^r(tn)<J_{2m+1}^r(tn)]\\
                         & =  &\lim_{r\rightarrow 0^+} \sum_{j=1}^{2m}\int_{\tilde{\Pi}_j^r}\mathbb{P}[J_{j-1}^r(tn)<J_{j}^r(tn)<J_{2m+1}^r(tn)|J_j^r=\pi]\mu_{J_j^r}(d\pi)\\
                         & = &\lim_{r\rightarrow 0^+} \sum_{j=1}^{2m}\int_{\tilde{\Pi}_j^r}\mathbb{P}[J_{j-1}^r(tn)<J_{j}^r(tn)|J_j^r=\pi] \times \\
                         & \ &  \mathbb{P}[J_{j}^r(tn)<J_{2m+1}^r(tn)|J_j^r=\pi]\mu_{J_j^r}(d\pi).
\end{eqnarray*}
Observe that whenever $J_j^r=\pi$, which occurs with positive probability, the events $[J_{j-1}^r(tn)<\pi(tn)]$ and $[\pi(tn)<J_{2m+1}^r(tn)]$ are independent. Indeed, conditionally on the event $[J_j^r=\pi]$, the occurrence of any event involving $J_{j-1}^r(tn)$ does not depend on the points of the Poisson processes to the right of $\pi$. In the same way, the occurrence of any event involving $J_{2m+1}^r(tn)$ does not depend on the points of the Poisson processes to the left of $\pi$. Therefore the last equality follows from the independence of Poisson processes on disjoint sets.

From Proposition \ref{FKG} we know that $\mathbb{P}[J_{j-1}^r(tn)<J_{j}^r(tn)|J_j^r=\pi]$ is increasing in $\pi$ and that
$\mathbb{P}[J_{j}^r(tn)<J_{2m+1}^r(tn)|J_j^r=\pi]$ is decreasing in $\pi$. Thus, from the $FKG$ inequality for $\mu_{J_j^r}(d\pi)$ we may
conclude that

\begin{eqnarray*}
\mathbb{P}[\eta_0^{\prime}\geq 3]& \leq & \lim_{r\rightarrow 0^+}\sum_{j=1}^{2m}\int_{\tilde{\Pi}_j^r}\mathbb{P}[J_{j-1}^r(tn)<J_{j}^r(tn)|J_j^r=\pi]\mu_{J_j^r}(d\pi)\\
&&\times\int_{\tilde{\Pi}_j^r}\mathbb{P}[J_{j}^r(tn)<J_{2m+1}^r(tn)|J_j^r=\pi]\mu_{J_j^r}(d\pi)\\
  &=& \lim_{r\rightarrow 0^+}\sum_{j=1}^{2m}\mathbb{P}[J_{j-1}^r(tn)<J_{j}^r(tn)]\mathbb{P}[J_{j}^r(tn)<J_{2m+1}^r(tn)]\\
  & = & \sum_{j=1}^{2m}\mathbb{P}[{\hat{S}^n}_{j-1}(tn)<{\hat{S}^n}_{j}(tn)]\mathbb{P}[{\hat{S}^n}_{j}(tn)<{\hat{S}^n}_{2m+1}(tn)]\\
  & \leq & \sum_{j=1}^{2m}\mathbb{P}[\hat{S}^n_{j-1}(tn)<\hat{S}^n_{j}(tn)]\mathbb{P}[\hat{S}^n_{0}(tn)<\hat{S}^n_{2m+1}(tn)]\\
  & \leq & 2m\mathbb{P}[\hat{S}^n_{(0,0)}(tn)<\hat{S}^n_{(1,0)}(tn)]\mathbb{P}[\hat{S}^n_{0}(tn)<\hat{S}^n_{2m+1}(tn)]\\
  & \leq & 4\frac{e\epsilon}{\alpha}\sqrt{n}\mathbb{P}[\hat{S}^n_{(0,0)}(tn)<\hat{S}^n_{(2/\alpha,0)}(tn)]\mathbb{P}[\hat{S}^n_{0}(tn)<\hat{S}^n_{2m+1}(tn)].
\end{eqnarray*}
If $\hat{S}^n_{(0,0)}(tn)<\hat{S}^n_{(2/\alpha,0)}(tn)$, then the paths starting at 
a distance $2/\alpha$ apart at level $k=0$ have not coalesced after having visited $C(\alpha)^{-1}tn$ levels. Then, by the Lemma
\ref{tempo2}

\begin{equation*}
\mathbb{P}[\hat{S}^n_{(0,0)}(tn)<\hat{S}^n_{(2/\alpha,0)}(tn)]\leq C\frac{\left\lceil 2/\alpha\right\rceil}{\sqrt{C(\alpha)^{-1}nt}}.
\end{equation*}
Analogously,

\begin{equation*}
\mathbb{P}[\hat{S}^n_{0}(tn)<\hat{S}^n_{2m+1}(tn)]\leq C\frac{\left\lceil \epsilon\sqrt{n} \right\rceil}{\sqrt{C(\alpha)^{-1}nt}}.
\end{equation*}
Therefore, $\mathbb{P}[\eta_0^{\prime}\geq 3]\leq \tilde{C}(\alpha,t)\epsilon^{2}$, where 
$\tilde{C}(\alpha,t)>0$. From (\ref{casi1}), we get

\begin{equation*}
\mathbb{P}[\eta_{S^n}(0,nt;0,\sqrt{n}\epsilon)\geq 3]\leq \tilde{C}(\alpha,t)\epsilon^{2}+\mathbb{P}[A_0^c].
\end{equation*}
Analogously, we can show for $k=1,2,...,r_n(t)$ that

\begin{equation*}
\mathbb{P}[\eta_{S^n}(nm_k,nt;0,\sqrt{n}\epsilon)\geq 3]\leq \tilde{C}(\alpha,t)\epsilon^{2}+\mathbb{P}[A_k^c].
\end{equation*}
Therefore,

\begin{equation*}
\sup_{0\leq k\leq r_n(t)}\mathbb{P}[\eta_{S^n}(nm_k,nt;0,\sqrt{n}\epsilon)\geq 3]\leq 
\tilde{C}(\alpha,t)\epsilon^{2}+\sup_{0\leq k\leq r_n(t)}\mathbb{P}[A_k^c].
\end{equation*}
Then, it follows from Lemma \ref{cansado} that

\begin{equation}\label{jojojo}
\limsup_{n\rightarrow\infty}\sup_{0\leq k\leq r_n(t)}
\mathbb{P}[\eta_{S^n}(nm_k,nt;0,\sqrt{n}\epsilon)\geq 3]\leq \tilde{C}(\alpha,t)\epsilon^{2}.
\end{equation}
Finally, condition $B_2$ follows from (\ref{jojojo}).

\section{Weak Convergence to Coalescing Brownian Motions}

\label{sec:I}
This section is devoted to the verification of condition $I$ stated in Theorem \ref{omega4}. The main ideas behind the proof of condition $I$ are essentially the same as those used in \cite{CFD}. However, some technical
difficulties arise in our model. Namely, the jump times are not homogeneous, the paths are restricted to a rectangle and the distribution of the increments are
independent but not identically distributed since their distribution are level dependent.

\begin{thm}
\label{condicionI} Let $(y_0,s_0), (y_1,s_1),\ldots ,(y_k,s_k)$ be $k+1$ different points in $\mathbb{R}\times[0,1/\alpha -1)$ such that $s_0\leq s_1\leq...\leq s_k $
and if $s_{i-1}=s_i$ for some $i$, $i=1,...,k$, then $y_{i-1}<y_i$. We have

\begin{equation}
\{ \beta^n_{(y_i^n,t_i^n)}, i=0,...,k  \}
\overset{\mathcal{D}}{\underset{n\rightarrow \infty}{\longrightarrow}}\{W^{(i)}, i=0,1,...,k \},
\end{equation}
where $(y_i^n,t_i^n)\in S^{(n)}$ for $i=1,2,\ldots,k$ and $W^{(i)}$ are coalescing Brownian motions with positive diffusive constant 
$2$ , starting at $\{(y_0,s_0), (y_1,s_1), \ldots, (y_k,s_k)\}$ restricted to the set $\mathbb{R}\times[0,1/\alpha-1]$.
\end{thm}

\noindent {\bf{Proof of Theorem \ref{condicionI}}}. The proof falls naturally into three steps.

\subsubsection*{Step I: Convergence of a single path}
  
From the proof of this result it will become apparent why the map used in the statement of Theorem \ref{omega} is needed. In virtue of the alternative model introduced in section \ref{sec:auxresults} it suffices to show convergence of a single path in $\beta^n$. Indeed, convergence of a single path follows immediately from the following

\begin{thm} \label{convuna} 
Let $(x_0^{\prime},t_0^{\prime})\in\mathbb{R}\times[0,1/\alpha-1]$. There exists a sequence $(\beta^n_{(x_0^{n},t_0^{n})})$ such that

\begin{eqnarray*}
\beta^n_{(x_0^{n},t_0^{n})} & \overset{\mathcal{D}}{\underset{n\rightarrow\infty}{\longrightarrow}} & \left\{\left(  x_0^{\prime}+\sqrt{2}W\left[ t-t_0^{\prime} \right], t^{\prime}\right)\right\}_{t_0^{\prime}\leq t\leq \frac{1}{\alpha}-1},
\end{eqnarray*}
\end{thm}

\begin{proof}
It follows directly from Lemma \ref{path} and the fact that $T(\xi^n)\overset{D}{=}\beta^n$.
\end{proof}

\subsubsection*{Step II}

Here we consider the case in which $k\geq 1, s_0=s_1=...=s_k$ and $y_{i-1}<y_i$ for every $i=1,2,...,k$. We focus on the case $s_0=0$. The proof for any other $s$ in $[0,1/\alpha-1)$ is entirely analogous.

Since, for any $x\in\mathbb{R}$, the probability of the event that $[-\sqrt{n}+x,x+\sqrt{n}]$ has no points of a Poisson Process with constant rate decays exponentially fast, we can safely assume that the starting points of the trajectories under consideration are deterministic. When $k=1$, in virtue of the (spatial) translation invariance of the Poisson Process , it is enough to consider the case of two trajectories starting from $(0,0)$ and $(y_1,0)$ for some $y_1 >0$ . Without loss of generality we assume that we are working in a probability space in which

\begin{equation}\label{im1}
\beta^n_{(y_1,0)}\overset{a.s}{\underset{n\rightarrow\infty}{\longrightarrow}}W^{(1)}.
\end{equation}
When $k=1$, the result follows from the following proposition, which is also a building block in the general case.

\begin{prop}\label{propo}
The conditional distribution of $\beta_{(0,0)}^{n}$
given $\beta_{(y_1,0)}^{n}$ converges a.s. to the distribution of $W^{(0)}$ given $W^{(1)}$.
\end{prop}

Let us write

\begin{equation}\label{im2}
K_n^{(i)}=\left\{ (\sqrt{n}x_1,nx_2): (x_1,x_2)\in \beta^n_{(y_i,0)}  \right\},\quad i=0,1
\end{equation}
and

\begin{equation}
K_n^{(i)}(l_j)=\sqrt{n}y_i+\sum_{h=1}^{j}S_h^{(i,n)},\hspace{.1cm} l_j=\frac{nj}{n-j}, \hspace{.1cm}1\leq j\leq {k_n}
\end{equation}
where $S_h^{i,n}$, $0\leq h\leq k_n+1$ are independent random variables such that $S_h^{(i,n)}\overset{\mathcal{D}}{=}{\arg\min}_{x\in Q_h}|x|$.

The idea of the proof in this step $II$ is inspired in an analogous result given in \cite{CFD}. However, some new ideas are needed in our case since we are working in the region $[0,n(1/\alpha-1)]$. The approach we have adopted in this work is the following:

Approximate the paths ${K_n}^{(0)}$ and ${K_n}^{(1)}$ before the time when they first come close (at a distance of order $n^{\gamma}$, $\gamma<1/2$) 
by independent paths. If the time when they first come close is less than $n(1/\alpha-1)-n^{\beta}$, with $\beta(\gamma)<1$ then they coalesce quickly. If not, the diffusive rescaling solves the problem.

Now, for $i=0,1$ let $\left\{ \tilde{S}_h^{(i,n)},0\leq h \leq k_n,i=0,1 \right\}$ be independent copies of random variables $S_h^{(i,n)}$. Define new random variables $\hat{S}_h^{(i,n)}$ by

\[
\hat{S}_h^{(i,n)}=
\begin{cases}
S_h^{(i,n)}, & \text{if $|S_h^{(i,n)}|\leq n^{\gamma}$}\\
\tilde{S}^{(i,n)}, &\text{otherwise}.
\end{cases}
\]
For each $i=0,1$, define a random path $\hat{K}_n^{(i)}$ by linear interpolation of the points

\begin{equation}\label{hat}
\hat{K}_n^{(i)}(l_k)=\sqrt{n}y_i+\sum_{h=1}^{k}\hat{S}_h^{(i,n)}
\end{equation}
for $k=1,2,...,k_n+1$, restricted to the region $[0,n(1/\alpha-1)]$. Then, define

\begin{equation}
\hat{\tau}_n=\min\{ l_k : \hat{K}_n^{(1)}(l_k)-\hat{K}_n^{(0)}(l_k)\leq 3n^{\gamma} \}\wedge k_n \nonumber
\end{equation}
and

\begin{equation}
{\tau}_n=\min\{ l_k : {K}_n^{(1)}(l_k)-{K}_n^{(0)}(l_k)\leq 3n^{\gamma} \}\wedge k_n. \nonumber
\end{equation}
Also, let $\tilde{K}_n^{(i)}$, $i=0,1$ be the random path defined by linear interpolation of the points

\[
\tilde{K}_n^{(i)}(l_k)=
\begin{cases}
\hat{K}_n^{(i)}(l_k), & \text{if $(l_k\leq \hat{\tau}_n)$}\\
\hat{K}_n^{(i)}(\hat{\tau}_n)+\sum_{h=g(\hat{\tau}_n)+1}^{k}\breve{S}_h^{(i,n)}, & \text{in other case}\\
\end{cases}
\]
for $k=1,2,...,k_n+1$, restricted to the region $[0,n(1/\alpha-1)]$ where $g(\hat{\tau}_n)=\frac{n\hat{\tau}_n}{n+\hat{\tau}_n}$ and

\[
\breve{S}_h^{(i,n)}=
\begin{cases}
\tilde{S}_h^{(0,n)}, &\text{if $i=0$}\\
S_h^{(0,n)}, & \text{if $i=1$.}\\
\end{cases}
\]

\begin{lem} \label{primerlema} 
Let $A_n:=\left\{ S_h^{(i,n)}=\hat{S}_h^{(i,n)}, h=1,2,...,k_n, i=0,1  \right\}$. Then,

\begin{equation}
\mathbb{P}[A_n^c \ \mbox{{\it i.o}} \ ]=0. \nonumber
\end{equation}
\end{lem}

\begin{proof}
Since

\begin{eqnarray}
\mathbb{P}[A_n^c] &\leq&  \sum_{i=0}^1\sum_{k=1}^{k_n}\mathbb{P}[S_k^{(i,n)}\neq \hat{S}_k^{(i,n)}] \nonumber \\
&=&   2\sum_{k=1}^{k_n}\mathbb{P}[|S_k^{(1,n)}|>n^\gamma] \nonumber \\
% =    2\sum_{k=1}^{k_n}e^{-\frac{n-k}{n}n^{\gamma}}
&\leq& 2ne^{-\alpha n^{\gamma}}, \nonumber
\end{eqnarray}
the result follows from the Borel-Cantelli lemma.
\end{proof}

\begin{cor} \label{borell} 
For $n$ large enough, let 
\[
B_n:=\left\{  K_n^{(i)}(l_k)=\hat{K}_n^{(i)}(l_k); k=1,2,...,k_n;i=0,1\right\}.
\]
Then,
\begin{equation}
\mathbb{P}[B_n^c \ \mbox{{\it i.o}} \ ]=0.
\end{equation}
\end{cor}

\begin{proof}
 The proof is immediate from the previous lemma and definition \ref{hat}.
\end{proof}

\begin{rem} \label{Observation} 
It follows from the construction of the processes $\tilde{K}_n^{(0)}$ and $\tilde{K}_n^{(1)}$  and the Borel-Cantelli lemma that, for $n$ large enough: (i) they are independent; and (ii)  $\tilde{K}_n^{(i)}\overset{\mathcal{D}}{=}K_n^{(i)}$, $i=0,1$.
\end{rem}
Let $\beta\in(2\gamma,1)$ and $l_n^{\beta}:=\min\{l_k : l_k\geq l_{k_n}-n^{\beta}\}$.

\begin{lem}
\label{taus} $\mathbb{P}[K_n^{(1)}(s)>K_n^{(0)}(s), \tau_n\leq s\leq \tau_n+n^{\beta}, \{\tau_n\leq l_n^{\beta}  \}   ]\underset{n\rightarrow\infty}{\longrightarrow}0$.
\end{lem}

\begin{proof}
Let $k_n^{\beta}\in\mathbb{N}$ be such that $l_n^{\beta}=\frac{nk_n^{\beta}}{n-k_n^{\beta}}$. Then,

\begin{multline}
\mathbb{P}\left[K_n^{(1)}(s)> K_n^{(0)}(s),\tau_n\leq s\leq \tau_n+n^{\beta}, \{\tau_n\leq l_n^{\beta}\}\right] \nonumber \\
%= \sum_{j=1}^{k_n^{\beta}}\mathbb{P}\left[K_n^{(1)}(s)> K_n^{(0)}(s),\tau_n\leq s\leq \tau_n+n^{\beta} ; \tau_n=l_j\right] \nonumber \\
\leq \sum_{j=1}^{k_n^{\beta}}\mathbb{P}\left[K_n^{(1)}(s)> K_n^{(0)}(s),\tau_n\leq s\leq \tau_n+n^{\beta}|\tau_n=l_j\right] \times \mathbb{P}[\tau_n=l_j]. \nonumber
\end{multline}
A straightforward computation guarantees the existence of a constant $C(\alpha)>1$ such that for all $k, k=0,1,...,k_n$

\begin{equation}\label{niveles}
1\leq l_{k+1}-l_k\leq C(\alpha).
\end{equation}
Moreover, $\mathbb{P}\left[K_n^{(1)}(s)> K_n^{(0)}(s),\tau_n\leq s\leq \tau_n+n^{\beta}\vert \tau_n=l_j\right]$ is bounded above by the probability of the following
event: there exist two fixed points at distance $\left\lceil 3n^{\gamma}\right\rceil$ at level $l_j$ such that the corresponding paths do not coalesce for a period of time of order $n^{\beta}$. From (\ref{niveles}) we conclude that the paths must cross along, at least, $\left\lfloor C(\alpha)^{-1}n^\beta\right\rfloor$ levels without coalescing. Analogously to Lemma \ref{tempo2} we have

\[
\mathbb{P}\left[\left.K_n^{(1)}(s)> K_n^{(0)}(s),\tau_n\leq s\leq \tau_n+n^{\beta}\right|\tau_n=l_j\right]\leq C\frac{\left\lceil 3n^{\gamma}\right\rceil}{\sqrt{\left\lfloor C(\alpha)^{-1}n^\beta \right\rfloor}}.
\]
We stress that this bound is valid for every $j=1,2,...,k_n^{\beta}$. Then,

\[
\mathbb{P}\left[K_n^{(1)}(s)> K_n^{(0)}(s),\tau_n\leq s\leq \tau_n+n^{\beta}, \{\tau_n\leq l_n^{\beta}\}\right]\leq C\frac{\left\lceil 3n^{\gamma}\right\rceil}{\sqrt{\left\lfloor C(\alpha)^{-1}n^\beta \right\rfloor}}.
\]
Since $\beta\in(2\gamma,1)$, the proof is complete. 
\end{proof}

\begin{cor} \label{clave}
Let $\tilde{Z}_n^{(i)}:=\left\{\left(\frac{x_1}{\sqrt{n}} ,\frac{x_2}{n}\right): (x_1,x_2)\in\tilde{K}_n^{(i)}     \right\}$. Then,

\begin{enumerate}
	\item $\tilde{Z}_n^{(1)}\overset{q.c}{\underset{n\rightarrow\infty}{\longrightarrow}}W^{(1)}$;
	\item $(\tilde{Z}_n^{(0)},\tilde{Z}_n^{(1)})\overset{\mathcal{D}}{\underset{n\rightarrow\infty}{\longrightarrow}}(\tilde{W}^{(0)},W^{(1)})$ where
	$\tilde{W}^{(0)}\overset{\mathcal{D}}{=}W^{(0)}$, and $\tilde{W}^{(0)}$ and $W^{(1)}$ are independent.
	\item The conditional distribution of $\left\{ \tilde{Z}_n^{(0)}(t), t\leq \frac{\hat{\tau}_n}{n}  \right\}$ given $\tilde{Z}_n^{(1)}$ converges
	a.s. to $\left\{ \tilde{W}^{(0)}(t), t\leq \tau  \right\}$ given $W^{(1)}$, where $\tau=\inf\{t\geq 0 : \tilde{W}^{(0)}(t)=W^{(1)}(t)\}\wedge\{1/\alpha-1 \}$.
\end{enumerate}
\end{cor}

\begin{proof} 

\begin{enumerate} 

\item 
The proof is immediate from (\ref{im1}), (\ref{im2}) and Corollary \ref{borell}.

\item 
The proof follows from part one of observation \ref{Observation} and Lemma \ref{convuna}.

\item 
Let $\bar{Z}_n^{(0)}$ be such that $\bar{Z}_n^{(0)}\overset{\mathcal{D}}{=}\tilde{Z}_n^{(1)}$ and

\begin{equation}
(\bar{Z}_n^{(0)},\tilde{Z}_n^{(1)})\overset{a.s}{\underset{n\rightarrow\infty}{\longrightarrow}}(\bar{W}^{(0)},W^{(1)}) \nonumber
\end{equation}
where $(\bar{W}^{(0)},W^{(1)})\overset{\mathcal{D}}{=}(\tilde{W}^{(0)},W^{(1)})$. Define,

\begin{equation}
\bar{\tau}_n:=\min\{l_k/n : \tilde{Z}_n^{(1)}(l_k/n)-\tilde{Z}_n^{(0)}(l_k/n)\leq 3n^{\gamma-1/2}   \}\wedge\{1/\alpha-1\} \nonumber
\end{equation}
and

\begin{equation}
\bar{\tau}:=\inf\{ t\geq 0: \bar{W}^{(0)}(t)=W^{(1)}(t)   \}\wedge\{1/\alpha-1\}. \nonumber
\end{equation}

It suffices to show that $\bar{\tau}_n\overset{a.s}{\underset{n\rightarrow\infty}{\longrightarrow}}\tau$. Then, the claim follows from

\begin{itemize}
	\item $\mathbb{P}\left[\inf\{ t\geq 0: \bar{W}^{(0)}=W^{(1)}(t)   \}=1/\alpha-1\right]=0$ (whose proof is immediate) and
	\item $\mathbb{P}\left[\left.\forall \epsilon>0, \bar{W}^{(0)}(t)>W^{(1)}(t) \hspace{.1cm}\text{for some $t\leq \bar{\tau}+\epsilon$}\right| \bar{\tau}<1/\alpha-1\right]=1$, \newline whose proof follows from the strong Markov property and the independence of the involved processes.
\end{itemize}

and the following deterministic result whose proof is left as an exercise.

\begin{lem}. Let $f_n,g_n,f,g$ be continues functions from $[0,1/\alpha-1)$ into $\mathbb{R}$ such that$f(0)<g(0)$, and let $T=\inf\{ t\geq 0: f(t)=g(t)  \}\in (0, 1/\alpha-1)$ be finite. Suppose that the functions above have the property that, for every $\delta>0$, there exists $t\in [T,T+\delta]\cap [0,1/\alpha-1]$ with $f(t)>g(t)$. Also, suppose that

\[
\sup_{0\leq t\leq (T+1)\wedge \{1/\alpha-1\}}|f_n(t)-f(t)|\underset{n\rightarrow\infty}{\longrightarrow}0
\]
and 
\[
\sup_{0\leq t\leq (T+1)\wedge \{1/\alpha-1\}}|g_n(t)-g(t)|\underset{n\rightarrow\infty}{\longrightarrow}0.
\]
Let $T_n=\inf\{t\geq 0: g_n(t)-f_n(t)\leq a_n \}$ where $a_n\geq 0$ is a given sequence such that $a_n\underset{n\rightarrow\infty}{\longrightarrow}0$. Then, $T_n\underset{n\rightarrow\infty}{\longrightarrow}T$.

\end{lem}
\end{enumerate}
\end{proof}

Now, define $Z_n^{(i)}:=\left\{\left(\frac{x_1}{\sqrt{n}},\frac{x_2}{n}   \right): (x_1,x_2)\in K_n^{(i)}   \right\}=\beta_{(y_i,0)}^{n}$, $i=0,1$.

\begin{cor} \label{IME} The conditional distribution of $\left\{{Z}_n^{(0)}(t), t\leq \frac{{\tau}_n}{n} \right\}$ given ${Z}_n^{(1)}$ converges a.s. to that of $\{\tilde{W}^{(0)}, t\leq \tau \}$ given $W^{(1)}$.   
\end{cor}

\begin{proof}
The result follows from Corollary \ref{borell} and part $(3)$ of Corollary \ref{clave}.
\end{proof}

\begin{cor}\label{IME1}
The process $\left\{Z_n^{(0)}, t\geq \frac{\tau_n}{n}\right\}$ converges in probability to $\{W^{(1)}(t), t\geq \tau \}$. 
\end{cor}

\begin{proof}
It suffices to to show that for any $\epsilon > 0$,

\begin{equation} \label{Cafe}
\mathbb{P}\left[ \max_{0\leq l\leq n^{\beta}}(K_n^{(1)}(l+\tau_n)-K_n^{(0)}(l+\tau_n))> 
\epsilon\sqrt{n} \right]\underset{n\rightarrow \infty}{\longrightarrow}0,
\end{equation}
where $\max_{0\leq l\leq n^{\beta}}(K_n^{(1)}(l+\tau_n)-K_n^{(0)}(l+\tau_n))=0$ if $l+\tau_n> k_n$. Now, conditioned  on $\tau_n=l_j$, the event

$$
\left\{\max_{0\leq l\leq n^{\beta}} (K_n^{(1)}(l+\tau_n)-K_n^{(0)}(l+\tau_n))> \epsilon\sqrt{n}\right\}
$$
is stochastically dominated by the increments of a process starting at level $l_j$ whose initial points are at a distance $3n^{\alpha}$ apart from each other. Then, by Doob's inequality for non-negative Martingales, we get

\[
\mathbb{P}\left[\left.\max_{0\leq l\leq n^{\beta}}(K_n^{(1)}(l+\tau_n)-K_n^{(0)}(l+\tau_n))
> \epsilon\sqrt{n}\right| \tau_n=l_j \right]\leq \frac{3n^{\gamma}}{\sqrt{\epsilon n}}.
\]
Then, (\ref{Cafe}) follows from the inequality above.
\end{proof}

\noindent {\bf Proof of Proposition \ref{propo}}. It follows from corollaries \ref{IME} and \ref{IME1}.

\vspace{.2cm}

Now we focus on the case $k>1$. It suffices to show that

\begin{equation}\label{final1}
\mathbb{E}[f_0(Z_n^{(0)})...f_k(Z_n^{(k)})]\underset{n\rightarrow\infty}
{\longrightarrow}\mathbb{E}[f_0(W^{(0)})...f_k(W^{(k)})]
\end{equation}
for every $f_0,...,f_k\in \mathcal{C}_b(\Pi_0^{1/\alpha-1},\mathbb{R}),$ 
the space of all real valued and bounded continues functions defined on $\Pi_0^{1/\alpha-1}$. Then,

\begin{eqnarray*}
\mathbb{E}[f_0(Z_n^{(0)})...f_k(Z_n^{(k)})]& = 
& \mathbb{E}\left\{ f_1(Z_n^{(1)})\mathbb{E}[f_0(Z_n^{(0)})f_2(Z_n^{(2)})...f_k(Z_n^{(k)})|Z_n^{(1)}] \right\}\\
                                           & = 
                                         & \mathbb{E}\left\{ f_1(Z_n^{(1)})
                                         \mathbb{E}[f_0(Z_n^{(0)})|Z_n^{(1)}]\mathbb{E}[f_2(Z_n^{(2)})...
                                         f_k(Z_n^{(k)})|Z_n^{(1)}] \right\}
\end{eqnarray*}
The second equality follows from the fact that, conditioned on $Z_n^{(1)}$, $Z_n^{(0)}$ and $(Z_n^{(2)}, \newline ...,Z_n^{(k)})$ are independent. Since $Z_n^{(1)}\overset{a.s}{\underset{n\rightarrow\infty}{\longrightarrow}}W^{(1)}$, the inductive hypothesis implies that
\begin{equation}
\mathbb{E}[f_0(Z_n^{(0)})|Z_n^{(1)}]\overset{q.c}{\underset{n\rightarrow\infty}{\longrightarrow}}\mathbb{E}[f_0(W^{(0)})|W^{(1)}] \nonumber
\end{equation}
and 
\begin{equation}
\mathbb{E}[f_2(Z_n^{(2)})...f_k(Z_n^{(k)})|Z_n^{(1)}]\overset{a.s}{\underset{n\rightarrow\infty}{\longrightarrow}}
\mathbb{E}[f_2(W^{(2)})...f_k(W^{(k)})|W^{(1)}]. \nonumber
\end{equation}
Since the functions $f_i, i=0,1,...,k$ are bounded, the dominated convergence theorem implies that

\begin{eqnarray}
\mathbb{E}\left\{ f_1(Z_n^{(1)})\mathbb{E}[f_0(Z_n^{(0)})|Z_n^{(1)}]\mathbb{E}[f_2(Z_n^{(2)})...f_k(Z_n^{(k)})|Z_n^{(1)}] \right\} \nonumber
\end{eqnarray}
converges, when $n$ goes to infinity, to
\begin{equation} \label{c}
%&\underset{n\rightarrow\infty} {\longrightarrow}&
\mathbb{E}\left\{ f_1(W^{(1)})\mathbb{E}[f_0(W^{(0)})|W^{(1)}] \mathbb{E}[f_2(W^{(2)})...f_k(W^{(k)})|W^{(1)}] \right\}.
\end{equation}
Note that (\ref{c}) equals

\begin{align}
\mathbb{E}\left\{ f_1(W^{(1)})\mathbb{E}[f_0(W^{(0)})f_2(W^{(2)})...f_k(W^{(k)})|W^{(1)}]\right\} = \\
\mathbb{E}[ f_1(W^{(1)})f_0(W^{(0)})f_2(W^{(2)})...f_k(W^{(k)})]. \nonumber
\end{align}
Then, (\ref{final1}) follows from the last convergence above. This completes the proof.

\subsubsection*{Step III: the general case}

The proof for the general case as stated in Theorem \ref{condicionI} is somewhat standard and follows in the same lines as those in the proof of Theorem $4$ of \cite{CFD}, page 1196 using step II proved above and the Markov property. For further details, we refer the reader to \cite{CFD}. 

\section{Verification of condition $B_1$}

We begin by stating that verifying condition $B_1$ of Theorem \ref{omega4} for the discrete radial Poissonian web is equivalent to show that, for every $t\in (0,1/\alpha-1)$,

\begin{equation}\label{fe}
\limsup_{n\rightarrow\infty}\sup_{(a,t_0)\in \mathbb{R}\times[0,1/\alpha-1]}\mathbb{P}[\eta_{\beta^n}(t_0,t;a,a+\epsilon)\geq 2]\underset{\epsilon\rightarrow 0_+}{\longrightarrow}0.
\end{equation}
For $t\in(0,1/\alpha-1)$ and $\epsilon>0$ fixed let $r_n(t):=\max\left\{k: \frac{k}{n-k}< 1/\alpha-1-t \right\}$ and $m_k:=\frac{k}{n-k}$. By the spatial translation invariance of the model we can safely remove the sup in (\ref{fe}). Then, (\ref{fe})  becomes

\begin{equation}\label{fe1}
\limsup_{n\rightarrow\infty}\sup_{0\leq k \leq r_n(t) }\mathbb{P}[\eta_{\beta^n}(m_k,t;0,\epsilon)\geq 2]\underset{\epsilon\rightarrow 0_+}{\longrightarrow}0.
\end{equation}

Let $S^n:=\left\{ (\sqrt{n}x_1,nx_2): (x_1,x_2)\in \beta^n  \right\}$. Then, for $k=0,...,r_n(t)$ we have

\[
\mathbb{P}[\eta_{\beta^n}(m_k,t;0,\epsilon)\geq 2]=\mathbb{P}[\eta_{S^n}(l_k,nt;0,\sqrt{n}\epsilon)\geq 2].
\]
Since the trajectories do not cross each other, a standard argument and Lemma \ref{tempo2} give, for $k=0,...,r_n(t)$,

\[
\mathbb{P}[\eta_{S^n}(l_k,nt;0,\sqrt{n}\epsilon)\geq 2]\leq C\frac{\epsilon \sqrt{n}}{\sqrt{C(\alpha)^{-1}nt}}=C_1(\alpha)\frac{\epsilon}{\sqrt{t}}.
\]
This gives (\ref{fe1}) and the proof is complete.

\section{Appendix} \label{sec:App}

For the sake of completeness we list some standard results used in the previous sections.

\subsection{Weak convergence}
In the first part of the Appendix we prove Lemma \ref{path}. First, we introduce some notation. Let $\alpha\in(0,1)$ be fixed. 
For $t\in[0,1-\alpha]$, make $f(t)=(1-t)$. Denote by $B$ a standard Brownian motion starting from the origin. Now, consider the process

\begin{equation}
(B[t^{\prime}],t^{\prime}), \quad t^{\prime}\in\left[0,\frac{1}{\alpha}-1\right]. \nonumber
\end{equation}
Then,

\begin{equation}
(B[t^{\prime}],t^{\prime})\overset{\mathcal{D}}{=}\left(B\left[ \frac{1}{f(t)}-1\right], \frac{1}{f(t)}-1    \right), \quad 0\leq t\leq 1-\alpha . \nonumber
\end{equation}
Also, consider a continuous mapping $T_1:C_{[0,1/\alpha-1]}\rightarrow C_{[0,1-\alpha]}$ such that

\begin{equation}
\left(B\left[ \frac{1}{f(t)}-1\right],\frac{1}{f(t)}-1\right)\overset{T_1}{\longmapsto}\left(cf(t)B\left[\frac{1}{f(t)}-1 \right],t\right),\quad 0\leq t\leq 1-\alpha, \nonumber
\end{equation}
where $c$ is a positive constant. Let $g(t) := \frac{1}{f(t)}-1$ and let $X$ be the process defined by

\begin{equation}
X_t := cf(t)B[g(t)],\quad t\in [0,1-\alpha]. \nonumber
\end{equation}
Let $\{ P_j \}_{j=1}^{\infty}$ be a sequence of independent Poisson processes on $\mathbb{R}$ with rate $1$. Define $\omega_j := {\arg\min}_{x\in P_j}|x|$ and let
$k_n(\alpha)=k_n=\lfloor n(1-\alpha)\rfloor$, where $n \in \mathbb{N}$.

Make $Y_n(0)=0$ and for $k=1,2,...,k_n$, make $Y_n(k) := f_n(k)\sum_{j=1}^k\frac{\sqrt{n}\omega_j}{n-j}$, where $f_n(k)=\frac{n-k}{n}$. Let

\begin{equation}\label{cafe}
X^n_t=Y_n(\lfloor nt\rfloor), 0\leq t\leq 1-\alpha.
\end{equation}

\begin{thm} \label{procesoprincipal} Let $X^n$ be the process defined by (\ref{cafe}). Consider the stochastic process defined by

\begin{equation}
Z^n_t=Y_n(\left\lfloor nt \right\rfloor)+(nt-\left\lfloor nt \right\rfloor)(Y_n(\left\lfloor nt \right\rfloor + 1)-Y_n(\left\lfloor nt \right\rfloor)),\quad t\in[0,1-\alpha].
\end{equation}
Then, $Z^n$ is an element of $C_{[0,1-\alpha]}$ and

\begin{equation}
Z_n\overset{\mathcal{D}}{\underset{n\rightarrow\infty}{\longrightarrow}} X.
\end{equation}
\end{thm}

\begin{prop}\label{m1}
For $n=1,2, \ldots$ let $\{Y_n(k),k=0,1,2,...\}$ be a discrete time, real-valued process. For $\alpha_n>0$ let

\begin{equation}\label{xn}
X^n_t := Y_n(\lfloor \alpha_n t \rfloor).
\end{equation}
Also, for every $t\in[0,T]$, let

\begin{equation}
Z^n_t := Y_n(\lfloor \alpha_nt \rfloor)+(\alpha_nt-\lfloor  \alpha_nt \rfloor)(Y_n(\lfloor \alpha_nt \rfloor+1)-Y_n(\lfloor \alpha_nt \rfloor)). \nonumber
\end{equation}
Then, $X_n$ converges in distribution to $X$ if, and only if, $Z_n$ converges in distribution to $X$.

\begin{equation}
X_n\overset{\mathcal{D}}{\underset{n\rightarrow\infty}{\longrightarrow}}X \Leftrightarrow Z_n\overset{\mathcal{D}}{\underset{n\rightarrow\infty}{\longrightarrow}}X. \nonumber
\end{equation}
\end{prop}

\begin{proof}
See \cite{ST} page 149.
\end{proof}

Note that $X^n$ and $X$ are elements of $D_{[0,1-\alpha]}$. Now we can prove that $X^n$ converges in distribution to $X$.

\begin{lem}\label{importante}
$X^n\overset{\mathcal{D}}{\underset{n\rightarrow\infty}{\longrightarrow}}X$.
\end{lem}

\noindent \textbf{Proof of Theorem \ref{procesoprincipal}} It is a simple consequence of Proposition \ref{m1} and the previous lemma.

Let $T > 0$ be fixed and let $D=D[0,T]$ denotes the space of all right continuous function $x: [0,T] \rightarrow \mathbb{R}$ with left limits. Denote by $\mathfrak{D}$ the corresponding Skorohod topology. Let $\mathbb{P}$ be a probability measure defined on $(D,\mathfrak{D})$. Denote by $T_{\mathbb{P}}$ the event where the following property holds almost surely:  $t (\in [0,T]) \in T_{\mathbb{P}}$ if, and only if, 
the projection  $\pi_t$ is continuous. Observe that $0$ and $T \in T_{\mathbb{P}}$.

\begin{thm}\label{convergencia}
Assume that, for every $t_1, t_2, ..., t_k \in T_{\mu_X},$ 

\begin{equation}
(X_{t_1}^n,...,X_{t_k}^n)\overset{\mathcal{D}}{\underset{n\rightarrow\infty}{\longrightarrow}} (X_{t_1},...,X_{t_k}),
\end{equation}
that 

\begin{equation}
X_T-X_{T-\delta}\overset{\mathcal{D}}{\underset{\delta\rightarrow 0}{\longrightarrow}}0 \nonumber
\end{equation}
and that, for $r\leq s\leq t, n\geq 1, \beta\geq 0, \alpha>1/2$

\begin{equation}
\mathbb{E}[|X_s^n-X_r^n|^{2\beta}|X_t^n-X_s^n|^{2\beta}]\leq [F(t)-F(r)]^{2\alpha} \nonumber
\end{equation}
where $F$ is a non-decreasing, continuous function on $[0,T]$. Then, $X^n\overset{\mathcal{D}}{\underset{n\rightarrow\infty}{\longrightarrow}}X$.
\end{thm}

\begin{proof}
See \cite{B}, page 142. 
\end{proof}

\begin{lem}\label{j1}
Assume that

\begin{equation}
(X_{t_1}^n,...,X_{t_k}^n)\overset{\mathcal{D}}{\underset{n\rightarrow\infty}{\longrightarrow}}(X_{t_1},...,X_{t_k}) \nonumber
\end{equation}
for any finite sequence of different points $t_1,t_2,...,t_k\in [0,1-\alpha]$. Then, $X^n\overset{\mathcal{D}}{\underset{n\rightarrow\infty}{\longrightarrow}}X$.
\end{lem}

\begin{proof} 
By Theorem \ref{convergencia} it suffices to show that

\begin{equation}\label{X}
X_1-X_{1-\delta}\overset{\mathcal{D}}{\underset{\delta\rightarrow 0}{\longrightarrow}}0.
\end{equation}
and that if $F$ is a non-decreasing, continuous function on $[0,1-\alpha]$ and if $r\leq s \leq t$, $n\geq 1$, $\beta\geq 0$, $\gamma>1/2$, then
\begin{equation} \label{Y}
\mathbb{E}[|X_s^n-X_r^n|^{2\beta}|X_t^n-X_s^n|^{2\beta}]\leq [F(t)-F(r)]^{2\gamma}.
\end{equation}
Now, (\ref{X}) follows from the almost surely continuity of the Brownian motion and the fact that the process $X$ is defined as a continuous mapping of a standard Brownian motion on $C_{[0,1/\alpha-1]}$. 

For $t\in[0,1-\alpha]$, let $\tilde{X}_t^n := \sum_{j=1}^{\lfloor nt\rfloor}\frac{\sqrt{n}\omega_j}{n-j}$. Make $J_0^n=0$. Thus, $X_t^{n}=f_n(\lfloor nt\rfloor)\tilde{X}_t^n$. Then,

\begin{eqnarray}\label{X_t}
|X_s^n-X_r^n|^2& = & |f_n(\lfloor ns\rfloor)\tilde{X}_s^n-f_n(\lfloor nr\rfloor)\tilde{X}_r^n|^2\nonumber\\
               & = & |f_n(\lfloor ns\rfloor)\tilde{X}_s^n-f_n(\lfloor ns\rfloor)\tilde{X}_r^n+f_n(\lfloor ns\rfloor)\tilde{X}_r^n-f_n(\lfloor nr\rfloor)\tilde{X}_r^n|^2\nonumber\\
               & = & |f_n(\lfloor ns\rfloor)(\tilde{X}_s^n-\tilde{X}_r^n)+\tilde{X}_r^n(f_n(\lfloor ns\rfloor)-f_n(\lfloor nr\rfloor))|^2\nonumber\\
               & \leq & 2|f_n(\lfloor ns\rfloor)|^2|\tilde{X}_s^n-\tilde{X}_r^n|^2 + 2|\tilde{X}_r^n|^2|f_n(\lfloor ns\rfloor)-f_n(\lfloor nr\rfloor)|^2\nonumber\\
               & \leq & 2|\tilde{X}_s^n-\tilde{X}_r^n|^2+2\left(\frac{\lfloor ns \rfloor-\lfloor nr\rfloor}{n} \right)^2|\tilde{X}_r^n|^2\nonumber\\
               & \leq & 2|\tilde{X}_s^n-\tilde{X}_r^n|^2+2\left(\frac{\lfloor nt \rfloor-\lfloor nr\rfloor}{n} \right)^2|\tilde{X}_r^n|^2.\nonumber\\
\end{eqnarray}
Analogously, 

\begin{eqnarray}\label{X_s}
|X_t^n-X_s^n|^2 & \leq & 2|\tilde{X}_t^n-\tilde{X}_s^n|^2+2\left(\frac{\lfloor nt \rfloor-\lfloor nr\rfloor}{n} \right)^2|\tilde{X}_s^n|^2.
\end{eqnarray}
If $t-r\leq 1/n$, then $\lfloor ns \rfloor=\lfloor nt\rfloor$ or $\lfloor ns\rfloor=\lfloor nr\rfloor$. Thus,

\begin{equation}\label{q1}
\mathbb{E}[|X_s^n-X_r^n|^{2}|X_t^n-X_s^n|^{2}]=0.
\end{equation}
Then we may assume that $t-r>1/n$. From (\ref{X_t}) and (\ref{X_s}) we get

\begin{equation}\label{q2}
|X_s^n-X_r^n|^2||X_t^n-X_s^n|^2\leq (I)+(II)+(III)+(IV),
\end{equation}
where

\begin{eqnarray}
(I) &=& 4|\tilde{X}_s^n-\tilde{X}_r^n|^{2}|\tilde{X}_t^n-\tilde{X}_s^n|^{2}, \nonumber \\
(II) &=& 4\left(\frac{\lfloor nt \rfloor-\lfloor nr\rfloor}{n} \right)^2|\tilde{X}_s^n-\tilde{X}_r^n|^{2}|\tilde{X}_s^n|^2, \nonumber \\
(III) &=& 4\left(\frac{\lfloor nt \rfloor-\lfloor nr\rfloor}{n} \right)^2|\tilde{X}_t^n-\tilde{X}_s^n|^{2}|\tilde{X}_r^n|^2 \ \mbox{and} \ \nonumber \\
(IV)&=&4\left(\frac{\lfloor nt \rfloor-\lfloor nr\rfloor}{n} \right)^4|\tilde{X}_r^n|^2|\tilde{X}_s^n|^2. \nonumber
\end{eqnarray}
By the independence of the increments we have

\begin{eqnarray*}
\mathbb{E}[(I)]& = & 4\mathbb{E}[|\tilde{X}_s^n-\tilde{X}_r^n|^{2}]\mathbb{E}[|\tilde{X}_t^n-\tilde{X}_s^n|^{2}]\nonumber\\
               & \leq & 4\left\{ \sum_{j=\lfloor nr\rfloor+1}^{\lfloor ns\rfloor}\frac{n\mathbb{E}[\omega_1^2]}{[n-j]^2} \right\}\times\left\{ \sum_{j=\lfloor ns\rfloor+1}^{\lfloor nt\rfloor}\frac{n\mathbb{E}[\omega_1^2]}{[n-j]^2}  \right\}\nonumber\\
               & \leq & 4(\mathbb{E}[\omega_1^2])^2\times\sum_{j=\lfloor nr\rfloor+1}^{\lfloor ns\rfloor}\frac{n}{(n\alpha)^2}\times\sum_{j=\lfloor ns\rfloor+1}^{\lfloor nt\rfloor}\frac{n}{(n\alpha)^2}\nonumber\\
               & \leq & 4(\mathbb{E}[\omega_1^2])^2\alpha^{-4}\left(  \frac{\lfloor ns\rfloor-\lfloor nr\rfloor}{n} \right)\left(  \frac{\lfloor nt\rfloor-\lfloor ns\rfloor}{n} \right)\nonumber\\
               & \leq & 4(\mathbb{E}[\omega_1^2])^2\alpha^{-4}\left(  \frac{\lfloor nt\rfloor-\lfloor nr\rfloor}{n} \right)^2.\nonumber\\
\end{eqnarray*}
Since $t-r>1/n$, we get

\begin{equation}\label{cero}
\mathbb{E}[(I)]\leq 4(\mathbb{E}[\omega_1^2])^2\alpha^{-4}\times 4(t-r)^2=C_1(t-r)^2.
\end{equation}
By H\"{o}lder inequality and the fact that $t-r>1/n$ we have

\begin{eqnarray}\label{11}
\mathbb{E}[(II)] & = & 4\left( \frac{\lfloor nt\rfloor-\rfloor nr\rfloor}{n} \right)^2\mathbb{E}[|\tilde{X}^n_s-\tilde{X}_r^n|^2|\tilde{X}_s^n|^2]\nonumber\\
                 & \leq & 16(t-r)^2\left\{ \mathbb{E}[|\tilde{X}^n_s-\tilde{X}_r^n|^4]  \right\}^{1/2}\left\{ \mathbb{E}[|\tilde{X}_s^n|^4] \right\}^{1/2}.\nonumber\\
\end{eqnarray}
Also,

\begin{eqnarray}\label{12}
\mathbb{E}[|\tilde{X}_s^n|^4] & \leq & \sum_{j=1}^{\lfloor ns\rfloor}\frac{n^2\mathbb{E}[\omega_1^4]}{[n-j]^4}+\underset{1\leq i,j\leq \lfloor ns\rfloor }{\sum_{i\neq j}}\frac{n^2(\mathbb{E}[\omega_1^2])^2}{[n-j]^4}\nonumber\\
                              & \leq & \frac{n^3\mathbb{E}[\omega_1^4]}{(n\alpha)^4}+\frac{n^4(\mathbb{E}[\omega_1^2])^2}{(n\alpha)^4}\leq C_2, \nonumber\\
\end{eqnarray}
where $C_2>0$ is a positive constant which depends on $\alpha, \mathbb{E}[\omega_1^2]$ and $\mathbb{E}[\omega_1^4]$. Analogously,

\begin{equation}\label{13}
\mathbb{E}[|\tilde{X}_s^n-\tilde{X}_r^n|^4]\leq C_2.
\end{equation}
Collecting (\ref{11}), (\ref{12}) and (\ref{13}) we get

\begin{equation}\label{segunda}
\mathbb{E}[(II)]\leq C_3(t-r)^2.
\end{equation}
Since $t-r>1/n$, we have

\begin{equation}\label{tercera}
\mathbb{E}[(III)]\leq C_3(t-r)^2.
\end{equation}
Using H\"{o}lder inequality, the same arguments used above and the fact that $t-r<1/n$ we get

\begin{equation*}
\mathbb{E}[(IV)]\leq C_4(t-r)^4.
\end{equation*}
If $|t-r|\leq 1$, then

\begin{equation}\label{cuarto}
\mathbb{E}[(IV)]\leq C_4(t-r)^2.
\end{equation}
Collecting (\ref{q1}), (\ref{q2}), (\ref{cero}) and (\ref{segunda}), (\ref{tercera}), (\ref{cuarto}) we get

\[
\mathbb{E}[|X_s^n-X_r^n|^2|X_t^n-X_s^n|^2]\leq C_5(t-r)^2.
\]
Then, (\ref{Y}) follows by making $\beta=1$, $\gamma=1$ and $F(t)=\sqrt{C_5}t$. This completes the proof.

\end{proof}

\begin{lem}\label{t2}
Let $c=\sqrt{2\mathbb{E}[\omega_1^2]}$. Then, the finite-dimensional distribution of the process $X^n$ converges weakly to the finite-dimensional distribution of the process $X$.
\end{lem}

\begin{proof} Let $t\in(0,1-\alpha]$ be fixed. Then, $X^n_t=f_n(\left\lfloor nt \right\rfloor)\tilde{X}_t^n$, where

\[
f_n(\left\lfloor nt \right\rfloor)=\frac{n-\left\lfloor nt \right\rfloor}{n}\quad\text{and}\quad \tilde{X}_t^n=\sum_{j=1}^{\left\lfloor nt \right\rfloor}\frac{\sqrt{n}\omega_j}{n-j}.
\]
For $k=1,2,...,k_n$, $r_n=\left\lfloor nt \right\rfloor$, let $Y_{nk} := \frac{\sqrt{n}\omega_k}{n-k}$. Observe that

\[
\tilde{X}_t^n=\sum_{j=1}^{r_n}Y_{nj}.
\]
Then,

\[
\mathbb{E}[Y_{nk}]=0, \quad \sigma_{nk}^2=Var[Y_{nk}]=\frac{n\mathbb{E}[(\omega_1)^2]}{[n-k]^2}\quad\text{and}\quad s_n^2=\sum_{j=1}^{r_n}\sigma_{nj}^2.
\]
Thus, $s_n^2=\sum_{j=1}^{r_n}\frac{n\mathbb{E}[(\omega_1)^2]}{[n-j]^2}\geq t\mathbb{E}[(\omega_1)^2]=C_0>0$. Let $\epsilon>0$ be fixed. Then,

\begin{eqnarray}\label{b1}
\sum_{j=1}^{r_n}\frac{1}{s_n^2}\mathbb{E}\left[Y_{nk}^21_{\{|Y_{nk}|\geq s_n\epsilon\}}  \right] & \leq & C_1\sum_{j=1}^{r_n}\mathbb{E}\left[\frac{n(\omega_j)^2}{[n-j]^2}1_{\left\{|\omega_j|\geq \frac{n-j}{\sqrt{n}}s_n\epsilon \right\}}\right]\quad (C_1=1/C_0) \nonumber\\
                & \leq & C_1\frac{n}{[n-r_n]^2}\sum_{j=1}^{r_n}\mathbb{E}\left[(\omega_j)^21_{\left\{|\omega_j|\geq \frac{n-j}{\sqrt{n}}s_n\epsilon \right\}}\right]\nonumber\\
                & \leq & C_1\frac{n}{[n-r_n]^2}\sum_{j=1}^{r_n}\mathbb{E}\left[(\omega_j)^21_{\left\{|\omega_j|\geq \frac{n-r_n}{\sqrt{n}}\sqrt{C_0}\epsilon \right\}}\right]\nonumber\\
                & = & C_1\frac{n}{[n-r_n]^2}r_n\mathbb{E}\left[(\omega_1)^21_{\left\{|\omega_1|\geq \frac{n-r_n}{\sqrt{n}}\sqrt{C_0}\epsilon . \right\}}\right]\nonumber\\
\end{eqnarray}
Then,

\begin{equation}\label{b2}
C_1\frac{n}{[n-r_n]^2}r_n=\frac{n\left\lfloor nt \right\rfloor}{t\mathbb{E}[(\omega_1)^2][n-\left\lfloor nt \right\rfloor]^2}\hspace{.2cm}\underset{n\rightarrow\infty}{\longrightarrow} \hspace{.2cm}\frac{1}{\mathbb{E}[(\omega_1)^2](1-t)},
\end{equation}

\begin{equation}\label{b3}
\frac{n-r_n}{\sqrt{n}}\sqrt{C_0}\epsilon \underset{n\rightarrow\infty}{\longrightarrow} +\infty \Rightarrow 1_{\left\{|\omega_1|\geq \frac{n-r_n}{\sqrt{n}}\sqrt{C_0}\epsilon \right\}}\overset{q.c}{\underset{n\rightarrow\infty}{\longrightarrow}} 0.
\end{equation}
It follows from the Dominated Convergence Theorem and (\ref{b1}), (\ref{b2}) and (\ref{b3}) that

\[
\sum_{j=1}^{r_n}\frac{1}{s_n^2}\mathbb{E}\left[Y_{nk}^21_{\{|Y_{nk}|\geq s_n\epsilon\}}  \right] 
\underset{n\rightarrow\infty}{\longrightarrow} 0.
\]
Therefore, the Lindeberg condition is satisfied. Thus, by Lindeberg's Theorem (see Billingsley \cite{Billingsley}, Theorem $27.2$, page 359) we have that

\begin{equation}\label{c_1}
\frac{\sum_{j=1}^{r_n}Y_{nj}}{s_n}\overset{\mathcal{D}}{\underset{n\rightarrow\infty}{\longrightarrow}} N(0,1).
\end{equation}
Note that $s_n^2=\sum_{j=1}^{\left\lfloor nt \right\rfloor}\frac{n\mathbb{E}[(\omega_1)^2]}{[n-j]^2}=2\mathbb{E}[\omega_1^2]\sum_{j=1}^{\left\lfloor nt \right\rfloor}\frac{1}{\left[1-\frac{j}{n}\right]^2}
\frac{1}{n}$. Therefore,

\begin{equation}\label{c_2}
s_n^2\underset{n\rightarrow\infty}{\longrightarrow} 2\mathbb{E}[\omega_1^2]\int_0^t\frac{1}{[(1-x)]^2}dx=2\mathbb{E}[\omega_1^2]g(t).
\end{equation}
Then, the result follows from Slutsky's Theorem, (\ref{c_1}) and (\ref{c_2}) that

\begin{equation*}
\sum_{j=1}^{r_n}Y_{nj}=s_n\frac{\sum_{j=1}^{r_n}Y_{nj}}{s_n}\overset{\mathcal{D}}{\underset{n\rightarrow\infty}{\longrightarrow}}\sqrt{2\mathbb{E}[\omega_1^2]} N(0,g(t)).
\end{equation*}
Thus,

\begin{equation}\label{d1}
\tilde{X}_t^n\overset{\mathcal{D}}{\underset{n\rightarrow\infty}{\longrightarrow}} \sqrt{2\mathbb{E}[\omega_1^2]}N(0,g(t)).
\end{equation}
Since $f_n(\left\lfloor nt \right\rfloor)\underset{n\rightarrow\infty}{\longrightarrow} (1-t)=f(t)$, \ref{d1}) and Slutsky's Theorem implies that

\begin{equation}\label{d2}
X_t^n=f_n(\left\lfloor nt \right\rfloor)\tilde{X}_t^n\overset{\mathcal{D}}{\underset{n\rightarrow\infty}{\longrightarrow}} cf(t) N(0,g(t)).
\end{equation}
where $c=\sqrt{2\mathbb{E}[\omega_1^2]}$.

Let $0<s<t\leq 1$ be fixed. We will show that

\begin{equation}
(X_s^n,X_t^n)\overset{\mathcal{D}}{\underset{n\rightarrow\infty}{\longrightarrow}}(cf(s)N(0,g(s)),cf(t)N(0,g(t))). \nonumber
\end{equation}
Consider the random vector

\begin{eqnarray}
({X}^n_s,{X}^n_t-{X}^n_s) & = & \left(f_n(\left\lfloor ns \right\rfloor)\tilde{X}_s^n,f_n(\left\lfloor nt \right\rfloor)\tilde{X}_t^n -f_n(\left\lfloor ns \right\rfloor)\tilde{X}_s^n \right)\nonumber\\
                          & = & \left(f_n(\left\lfloor ns \right\rfloor)\tilde{X}_s^n, f_n(\left\lfloor nt \right\rfloor)(\tilde{X}_t^n-\tilde{X}_s^n)\right) \nonumber \\
                          & \ & + \ (0,\tilde{X}_s^n(f_n(\left\lfloor nt \right\rfloor)-f_n(\left\lfloor ns \right\rfloor))).\nonumber
\end{eqnarray}
Analogously to the previous case we get

%\begin{eqnarray}
%\left(f_n(\left\lfloor ns \right\rfloor)\tilde{X}_s^n, f_n(\left\lfloor nt \right\rfloor)(\tilde{X}_t^n-\tilde{X}_s^n)\right)&\overset{\mathcal{D}}{\underset{n\rightarrow\infty}%{\longrightarrow}} &\left(cf(s)N_1(0,g(s), cf(t)N_2\left(0,\int_s^t\frac{1}{(1-x)^2}dx \right)  \right)\nonumber\\
%& = & \left(cf(s)N_1(0,g(s)), cf(t)N_2(0,g(t)-g(s))\right),\nonumber\\
%\end{eqnarray}

\begin{equation} \label{e1}
\left(f_n(\left\lfloor ns \right\rfloor)\tilde{X}_s^n, f_n(\left\lfloor nt \right\rfloor)(\tilde{X}_t^n-\tilde{X}_s^n)\right) \overset{\mathcal{D}}{\underset{n\rightarrow\infty}{\longrightarrow}} \left(cf(s)N_1(0,g(s)), cf(t)N_2(0,g(t)-g(s))\right),\\
\end{equation}
where we have used the independence of the coordinates of the random vector. Here $N_1$ and $N_2$ are independent and normal distributed random variables. We also have

\begin{eqnarray}\label{e2}
(0,\tilde{X}_s^n(f_n(\left\lfloor nt \right\rfloor)-f_n(\left\lfloor ns \right\rfloor))) & \overset{\mathcal{D}}{\underset{n\rightarrow\infty}{\longrightarrow}} & \left(0,(f(t)-f(s))c N_1(0,g(s))\right).
\end{eqnarray}
Therefore, from (\ref{e1}) and (\ref{e2}) we have that

\begin{align}\label{e4}
(X^n_s, X^n_t-X^n_s) \overset{\mathcal{D}}{\underset{n\rightarrow\infty}{\longrightarrow}} (cf(s)N_1(0,g(s)),cf(t)N_2(0,g(t)-g(s))) \nonumber \\ 
+(0,c[f(t)-f(s)]N_1(0,g(s))). 
\end{align}
Then, from (\ref{e4}) and the independence between $N_1$ and $N_2$, we get

%\begin{eqnarray}
%(X^n_s, X^n_t) & = & (X^n_s, X^n_t-X^n_s)+(0,X^n_s)\nonumber\\
%                               & \overset{\mathcal{D}}{\underset{n\rightarrow\infty}{\longrightarrow}} & (cf(s)N_1(0,g(s)),cf(t)N_2(0,g(t)-g(s)))+(0,c[f(t)-f(s)]N_1(0,g(s)))\nonumber\\
%                               &  &+ (0,cf(s)N_1(0,g(s)))\nonumber\\
%                               & = & (cf(s)N_1(0,g(s)),cf(t)N_2(0,g(t)-g(s)))+(0,cf(t)N_1(0,g(s)))\nonumber\\
%                               & = & (cf(s)N_1(0,g(s)), cf(t)N_3(0,g(t))).\nonumber\\
%\end{eqnarray}

\begin{equation}
(X^n_s, X^n_t) = (X^n_s, X^n_t-X^n_s)+(0,X^n_s) \overset{\mathcal{D}}{\underset{n\rightarrow\infty}{\longrightarrow}} (cf(s)N_1(0,g(s)), cf(t)N_3(0,g(t)))), \nonumber
\end{equation}
where $N_3\left(0, g(t)\right)$ is a zero mean normal random variable with variance $g(t)$. Let $m\in\mathbb{N}$ e $0\leq t_1\leq ...\leq t_m\leq 1$ be fixed. Then, analogously to the previous case, we get

\begin{equation}\label{e5}
(X_{t_1}^n,...,X_{t_m}^n)\overset{\mathcal{D}}{\underset{n\rightarrow\infty}{\longrightarrow}}(cf(t_1)N_1(0,g(t_1)),...,cf(t_m)N_m(0,g(t_m))).
\end{equation}
Now, the result follows easily from (\ref{e5}). 
\end{proof}

Now we prove Lemma \ref{importante}.

\noindent {\bf{Proof of Lemma \ref{importante}}} It follows immediately from Theorem \ref{convergencia} and Lemmas \ref{j1} and \ref{t2}.

We finish the appendix by giving the proof of Lemma \ref{path}.

\noindent {\bf{Proof of Lemma \ref{path}}} By Lemma \ref{1path} is enough to show that result for deterministic point $(-1,0)$. Now the lemma follows directly from Lemma \ref{importante}.

\end{document}